\documentclass[11pt]{article}
\usepackage[margin=1in]{geometry}

\usepackage{graphicx}
\usepackage{subcaption}
\usepackage{amsmath}
\usepackage{amssymb}
\usepackage{amsthm}
\usepackage{color}
\allowdisplaybreaks[1]

\usepackage{caption}
\usepackage{subcaption}
\usepackage{algorithm}
\usepackage{algpseudocode}
\usepackage{hyperref}

\theoremstyle{plain}
\newtheorem{theorem}{Theorem}[section]

\newtheorem{proposition}[theorem]{Proposition}

\usepackage{subcaption}

\theoremstyle{definition}

\theoremstyle{remark}

\def\0{{\bf 0}}
\def\1{{\bf 1}}

\def \Snn{ \mathbb{S}^{n \times n}}

\def \bmat{\left[\begin{matrix}}
	\def \emat{\end{matrix}\right]}
\def \bvec{\left(\begin{matrix}}
	\def \evec{\end{matrix}\right)}

\def \xy1vec{\left[\begin{matrix}x\\y\\1\end{matrix}\right]}
\def \QED{\begin{flushright}\Halmos\end{flushright}\end{proof}}
\def \defeq{\mathrel{\mathop{:}}=}
\def \Tr{\mathrm{Tr}}
\def \xbar{\bar{x}}
\def \ybar{\bar{y}}
\def \grad{\triangledown}

\def \gp3d{\grad p_3(d)}
\def \Hess{\triangledown^2}
\def \R{\mathbb{R}}
\def \Rn{\R^n}
\def \beq{\begin{equation}}
\def \eeq{\end{equation}}
\def \baeq{\begin{equation*}\begin{aligned}}
\def \eaeq{\end{aligned}\end{equation*}}
\newcommand{\baeql}[1]{\begin{equation}\label{#1}\begin{aligned}}
\def \eaeql{\end{aligned}\end{equation}}

\def \otn{\{1, \ldots, n\}}

\title{\LARGE \bf Third-Order Newton Methods with Semidefinite Programming}
\author{Jeffrey Zhang\thanks{UC Riverside, Department of Statistics. \texttt{jeffrey.zhang1@ucr.edu}}}
	
\begin{document}
\date{}
\maketitle
\begin{abstract}
	\noindent We study a third-order Newton method in which each iteration is obtained by minimizing the third-order Taylor expansion of the objective via a semidefinite programming (SDP) subproblem. In contrast to classical Newton methods, which rely on quadratic approximations, our approach leverages higher-order information to better capture the local geometry of the objective. We establish cubic convergence under standard assumptions, including strong convexity and sufficiently accurate initialization. We also propose practical modifications to address infeasibility of the SDP subproblem, including regularization and alternative objective formulations.
\end{abstract}


\section{Introduction}\label{Sec: Introduction}

Newton's method is one of the most fundamental algorithms in numerical optimization. With the goal of minimizing an objective function $f$, the premise is simple: given a current iterate, one approximates the objective function by its second-order Taylor expansion at that point and takes the global minimum of this approximation as the next iterate. A classical result states that, under certain assumptions—namely strong convexity of the objective function and sufficient proximity of the current iterate to the global minimum—this method achieves quadratic convergence. As a result, it serves as a core subroutine in many algorithms, including interior-point methods. However, these assumptions can be restrictive, and Newton's method may fail in various ways when they are not satisfied. Consequently, the literature on improving Newton's method is vast, often combining approaches such as Levenberg-Marquardt regularization~\cite{levenberg1944method,marquardt1963algorithm}, trust regions~\cite{Powell1970ANA, Mor1982RecentDI, conn2000trust}, damping~\cite{ortega2000iterative}, and cubic regularization~\cite{nesterov2006cubic}.

More recently, there has been interest in incorporating higher-order information into Newton-type methods. Higher order derivatives allow for approximations of the objective function, more accurately capturing it in a wider region. The primary bottleneck to utilizing higher order derivatives is the NP-hardness of minimizing quartic polynomials~\cite{murty1987some} (and onward), even locally~\cite{ahmadi2022complexityb}. Thus most of the work towards this have focused on using a convex regularization term to enforce convexity\cite{bubeck2019near,gasnikov2019optimal,nesterov2021implementable,cartis2023second,ahmadi2024higher}. One takes the $p$-th order Taylor expansion, and adds a term $\|\cdot\|^{p+1}$ multiplied by a sufficiently large scalar. The key observation is that with a sufficiently large regularization term, this function is convex. The iterates are then based on the minimum of each of these now convex functions, while they are high-degree polynomials, are amenable to convex optimization techniques. Note that finding this regularization scalar is nontrivial, as even deciding convexity is an NP-hard problem~\cite{ahmadi2013np}. The work in~\cite{ahmadi2024higher} gets around this by using the more tractable sum-of-squares convexity rather than convexity to produce a generalizable Newton's method that has polynomial-time iterations.

The focus of this paper will be on establishing the basic theory of unregularized third-order Newton methods, along with some regularization strategies. Our approach in this paper is based on leveraging the tractability of local minimization of cubic polynomials, rather than extending the convex-regularized approach. In particular, the local minimum of a cubic polynomial can be found by semidefinite programming (SDP)~\cite{ahmadi2022complexity}. Moreover, such an SDP can be written with only the coefficients of the cubic polynomial, in particular not including Lipschitz constants required for a regularization term.

For expanded analysis on use of this framework, we refer the reader to~\cite{cai2026globally}, which draws upon an earlier version of this manuscript. That work adapts this framework into an algorithm with global guarantees and provides additional convergence analysis and extensive simulation studies. The simulation studies we perform in this work will focus on different aspects compared to those in \cite{cai2026globally}.

Compared to convex-regularized third-order methods, our approach carries certain advantages. The first is implementability; the mapping from a problem and iterate to a subproblem is straightforward and doesn't require specially tuned parameters. For elaboration on this topic, see the discussion in Section 2.1 of \cite{cai2026globally}. The work in \cite{cai2026globally} found additional advantages empirically as well. When an unregularized step was possible, the lack of a regularization term allowed for better fidelity in the Taylor approximation, allowing for better navigation around complex curvature. The authors also found better robustness and reliability in low-dimensional nonconvex settings compared to their convex-regularized method.

Compared to existing second-order methods, the advantages of third-order methods do not lie in wall-clock time or per-iteration cost. The complexity-per-iteration will be significantly greater, as utilizing the third order derivatives involves first computing the derivatives themselves, and then solving SDPs rather than performing matrix inversions. Rather, the advantages of third-order methods lie in reduced iteration counts and improved convergence behavior (see Section 5 of \cite{cai2026globally} for supporting simulation studies). The usage of third-order methods should focus on scenarios where iteration complexity is of elevated importance, such as federated learning, where function evaluations are expensive, such as simulation studies, or where the objective function is highly nonconvex and second-order methods are not successful.

\subsection*{Organization and Contributions}\label{SSec: Organization}

We review preliminaries on local minima, cubic polynomials, and Taylor expansions in Section~\ref{Sec: Preliminaries}. In Section~\ref{Sec: Convergence}, we demonstrate the local well-definedness and cubic convergence of the unregularized third-order method (Theorem~\ref{Thm: Convergence}). In Section~\ref{Sec: Adaptations}, we propose strategies for when the algorithm is not well-defined. In Section~\ref{Sec: Simulations}, we provide numerical simulations illustrating the third-order Newton's method. We provide final commentary and directions for future work in Section~\ref{Sec: Conclusion}.

\section{Preliminaries and Notation}\label{Sec: Preliminaries}

Throughout this paper, we will be working primarily with cubic polynomials. For an $n$-variate cubic polynomial $p$, we will write it with the convention
\beq \label{Eq: Cubic Poly Form} p(x) = \frac{1}{6}\sum_{i=1}^n x^T x_iH_ix + \frac{1}{2}x^TQx + b^Tx + c,\eeq
where $b$ is an $n \times 1$ vector containing the linear coefficients of $p$, $Q$ is a symmetric $n \times n$ matrix whose $i,j$-th entry is $\frac{\partial^2 p}{\partial x_i\partial x_j}$, and each $H_i$ for $i \in \otn$  is an $n \times n$ matrix whose $j,k$-th entry is $\frac{\partial^3 p}{\partial x_i \partial x_j \partial x_k}$. In such notation, the gradient of the polynomial can be written as
\beq \label{Eq: Cubic Poly Form Gradient} \grad p(x) = \frac{1}{2}\sum_{i=1}^n x_iH_ix + Qx + b,\eeq
and its Hessian is
\beq \label{Eq: Cubic Poly Form Hessian} \grad^2 p(x) = \sum_{i=1}^n x_iH_i + Q.\eeq
Note that by convention, $Q$ is symmetric and the matrices $H_i$ satisfy
\begin{equation}\label{Eq: Valid Hessian}
(H_i)_{jk} = (H_j)_{ik} = (H_k)_{ij}.
\end{equation}
Using this, it is not difficult to show (see, e.g., \cite{ahmadi2022complexity}[Lemma 4.1]) that for two $n \times 1$ vectors $x$ and $y$, we have
\beq\label{Lem: Hessian switch}\sum_{i=1}^n x_iH_iy = \sum_{i=1}^n y_iH_ix = \bmat x^TH_1y \\ \vdots \\ x^TH_ny \emat.\eeq

A point $x$ is a \emph{(strict) local minimum} of a function $f$ if there exists an $\epsilon > 0$ for which $f(x) \le f(y)$ (resp. $f(x) < f(y)$) whenever $\|y - x\| < \epsilon$. For a continuously differentiable function, it is well-known that any such point $x$ must satisfy the \emph{first-order optimality condition} that the gradient $\grad f(x)$ is zero, and also the \emph{second-order optimality condition} that $\grad f(x) = 0$ and the Hessian $\Hess f(x)$ is positive semidefinite (psd), i.e. has nonnegative eigenvalues. Furthermore, if a point $x \in \Rn$ is a (strict) local minimum of a function $f: \Rn \to \R$, then for any direction $d \in \Rn$ (with $d \ne 0$), the restriction of $f$ in the direction $d$ ---i.e. the univariate function $q(\alpha) \defeq p(x + \alpha d)$---has a (strict) local minimum at $\alpha = 0$.

One can observe that a univariate cubic polynomial has either no local minima, exactly one local minimum (which is strict), or infinitely many non-strict local minima in the case that the polynomial is constant. A cubic polynomial $p(x) \defeq ax^3 + bx^2 + cx + d$ has a local minimum if and only if $b^2 - 3ac \ge 0$. If it exists, the local minimum is located at $$x = \frac{-b + \sqrt{b^2 - 3ac}}{3a} = -\frac{c}{b + \sqrt{b^2 - 3ac}},$$
and this can be deduced by finding the roots of the derivative. The inflection point is located at $x = -\frac{b}{3a}$, and the distance from the inflection point to the local minimum (and maximum), if they exist, is $|\frac{\sqrt{b^2 - 3ac}}{3a}|$. Note that in the case $a = 0$, the coordinate for the local minimum reduces to that of the quadratic case (and there is no inflection point).

Multivariate cubic polynomials also have either no local minima, exactly one strict local minimum, or have infinitely many non-strict local minima. This can be argued by observing that the restriction of a multivariate cubic polynomial is a univariate cubic polynomial, and that a (strict) local minimum of a function is also a (strict) local minimum along the restriction to any line. A local minimum of a cubic polynomial is strict if and only if the Hessian at that point is positive definite \cite[Theorem 3.1, Corollary 3.4]{ahmadi2022complexity}. If it exists, it can be found by solving an SDP. An SDP is an optimization problem of the form:

\begin{equation*}
\begin{aligned}
& \underset{X \in \mathbb{S}^{n \times n}}{\min}
& & \Tr(CX) \\
& \text{subject to}
& & \Tr(A_iX) = b_i, i = 1, \ldots, m\\
&&& X \succeq 0,
\end{aligned}
\end{equation*}
where $\Tr$ denotes the trace of a matrix, i.e., the sum of the diagonal entries, $\mathbb{S}^{n \times n}$ denotes the set of $n \times n$ symmetric matrices, $A \succeq B$ denotes that a matrix $A - B$ is positive semidefinite, and $C, A_i, i = 1,\ldots, n$ are matrix data, and $b_i$ are scalar data. It is well-known that semidefinite programs can be solved to arbitrary accuracy in polynomial time~\cite{vandenberghe1996semidefinite}. For a cubic polynomial written in the form (\ref{Eq: Cubic Poly Form}), the SDP that finds a local minimum can be given by

\begin{equation}\label{Eq: complete cubic SDP}
    \begin{aligned}
	& \underset{X \in \Snn, x \in \Rn, y \in \R}{\inf}
	& & 0\\
	& \text{subject to}
    && \frac{1}{2}\Tr(QX) + b^Tx + \frac{y}{2} = 0\\
	&&&  \frac{1}{2}\Tr(H_iX)+e_i^TQx+b_i=0, \forall i = 1, \ldots, n,\\
	&&& \bmat \sum_{i=1}^n x_iH_i + Q & \sum_{i=1}^n \Tr(H_iX)e_i+Qx \\ (\sum_{i=1}^n \Tr(H_iX)e_i+Qx)^T & y\emat \succeq 0,\\
	&&& \bmat X & x \\ x^T & 1 \emat \succeq 0,
    \end{aligned}
\end{equation}
where $e_i$ denotes the vector of zeros except the $i$-th entry, which is a one~\cite[Algorithm~2]{ahmadi2022complexity}. For practical implementation, one can minimize the expression $\frac{1}{2}\Tr(QX) + b^Tx + \frac{y}{2}$, as it is always nonnegative, and is zero if and only if the cubic polynomial has a second-order point~\cite{ahmadi2022complexity}.

For the implementations in this paper, we instead use the formulation
\begin{equation}\label{Eq: cubic SDP}
    \begin{aligned}
	& \underset{X \in \Snn, x \in \Rn, y \in \R}{\inf}
	& & \frac{1}{3}\Tr(QX) + b^Tx + \frac{y}{6}\\
	& \text{subject to}
	&&  \frac{1}{2}\Tr(H_iX)+e_i^TQx+b_i=0, \forall i = 1, \ldots, n,\\
	&&& \bmat \sum_{i=1}^n x_iH_i + Q & \sum_{i=1}^n \Tr(H_iX)e_i+Qx \\ (\sum_{i=1}^n \Tr(H_iX)e_i+Qx)^T & y\emat \succeq 0,\\
	&&& \bmat X & x \\ x^T & 1 \emat \succeq 0,
    \end{aligned}
\end{equation}
which only differs from (\ref{Eq: complete cubic SDP}) by the objective function. When the cubic polynomial has a local minimum, the two will share the same solution, as the quantity $\frac{1}{2}\Tr(QX) + b^Tx + \frac{y}{2}$ represents the duality gap between (\ref{Eq: cubic SDP}) and its dual. However, the formulation in (\ref{Eq: cubic SDP}) has greater interpretability, and also demonstrates slightly better performance empirically (see Section~\Ref{Sec: Simulations}). Both (\ref{Eq: complete cubic SDP}) and (\ref{Eq: cubic SDP}) can be derived as lifted SDP relaxations, of the following two optimization problems:

\begin{minipage}[t]{0.4\textwidth}
\begin{equation*}
    \begin{aligned}
	(\ref{Eq: complete cubic SDP}) \le & \underset{x \in \Rn}{\inf} 
	& & x^T\grad p(x)\\
	& \text{subject to}
	&&  \grad p(x) = 0, \Hess p(x) \succeq 0\\
    \end{aligned}
\end{equation*}
\end{minipage}\hfill
\begin{minipage}[t]{0.4\textwidth}
\begin{equation}\label{Eq: unrelaxed problems}
    \begin{aligned}
	(\ref{Eq: cubic SDP}) \le & \underset{x \in \Rn}{\inf} 
	& & p(x)\\
	& \text{subject to}
	&&  \grad p(x) = 0, \Hess p(x) \succeq 0\\
    \end{aligned}
\end{equation}
\end{minipage}\\

A \emph{norm} is a function defined on a vector space $\mathcal{X}$ satisfying
\begin{enumerate}
    \item $\|x\| \ge 0, \forall x \in \mathcal{X}$, with $\|x\| = 0$ if and only if $x = 0$ (positive definiteness),
    \item $\|\alpha x\| = |\alpha| \|x\|, \forall \alpha \in \R, x \in \mathcal{X}$ (linearity),
    \item $\|x + y\| \le \|x\| + \|y\|, \forall x,y \in \mathcal{X}$ (triangle inequality).
\end{enumerate}
A norm by extension obeys the ``reverse triangle inequality'', i.e. that $\|x - y\| \ge \|x\| - \|y\|$. The norms we use in this paper are standard symmetric $p$-linear norms. For vectors and matrices respectively, these are
$$\|x\| = \underset{\|y\| \le 1}{\max} x^Ty \text{\quad and \quad}\|A\| = \underset{\|x\|,\|y\| \le 1}{\max} x^TAy,$$
which reduce to the Euclidean norm and largest eigenvalue in absolute value. These norms are submultiplicative, that is, for vectors $x,y \in \Rn$ and matrices $A, B \in \R^{n \times n}$, we have $\|x^Ty\| \le \|x\|\|y\|, \|Ax\| \le \|A\|\|x\|$, and $\|AB\| \le \|A\|\|B\|$. Finally, for a symmetric $n \times n \times n$ tensor $H$ expressed as a set of matrices $\{H_1, \dots, H_n\} \subset \R^{n \times n}$ satisfying (\ref{Eq: Valid Hessian}), its norm is
\begin{equation}
    \|H\| = \underset{\|x\|, \|y\|, \|z\|\le 1}{\max} \sum_{i=1}^n x^Ty_iH_iz.
\end{equation}

By the Fundamental Theorem of Calculus, for a continuously differentiable function $F: \R^n \to \R^m$ with Jacobian $J$, we have
$$F(y) = F(x) + \int_0^1 J(x + \tau (y -x))(y -x)d\tau.$$
Integrals also obey the following inequality:
$$\|\int_x^y F(\tau)d\tau\| \le \int_x^y \|F(\tau)\|d\tau.$$

A function $f: \mathcal{R} \to \R$ is \emph{Lipschitz continuous} if there exists a constant $L$ such that for any $x, y \in \mathcal{R}$, $\|f(x) - f(y)\| \le L\|x-y\|$. A function $f: \Rn \to \R$ is called \emph{strongly convex} if there exists a constant $m > 0$ such that $\Hess f(x) \succeq mI, \forall x \in \Rn$, where $I$ denotes the identity matrix. For a given function $f$ and a point $\xbar$, we denote by $\Phi_{f,\xbar}$ the third-order Taylor expansion of $f$ around $\xbar$, i.e. the function
$$\sum_{i=1}^n (x-\xbar)^T (x-\xbar)_i \grad_i^3 p(\xbar) (x-\xbar) + \frac{1}{2} (x-\xbar)^T \grad^2 p(\xbar) (x-\xbar) + (x-\xbar)^T\grad p(\xbar).$$

\section{Convergence of the Unregularized Third-Order Newton's Method}\label{Sec: Convergence}

Following the basic premise for the classical quadratic Newton's method, consider the following algorithm for minimizing a function:

\begin{algorithm}[H]
	\caption{Unregularized Third-Order Newton's Method}\label{Alg: 3ON}
	\begin{algorithmic}[1]
		\State {\bf Input:} A function $f$, threshold $\epsilon > 0$
		\State Initialize a point $x^0$, $k = 0$
		\State {\bf While} $\|\grad f(x^k)\| > \epsilon$
		\State \quad Let $x^{k+1}$ be the local minimum of $\Phi_{f,x^k}$
		\State \quad Set $k = k+1$
        \State {\bf Output:} $x^k$
	\end{algorithmic}
\end{algorithm}

\noindent Note that Step 4 of Algorithm~\ref{Alg: 3ON} is solving an SDP. Our goal in this section is to prove a result that is an analogue to the well-known result for the quadratic Newton's Method, that is, if the function we are minimizing is strongly convex and we initialize close enough to a local minimum, then this algorithm converges cubically.

Unlike for quadratic polynomials, the Hessian being positive definite at the current iterate is not sufficient for a well-defined next iterate. We therefore begin by showing that this algorithm is well-defined, at least if the current iterate is sufficiently close to a local minimum.

\begin{theorem}\label{Thm: Well-defined}
Suppose that a strongly convex function $f$ has a strict local minimum $\xbar$. Then there exists $R > 0$ such that for all $x$ satisfying $\|x - \xbar\| < R$, $\Phi_{f,x}$ has a strict local minimum.
\end{theorem}

\begin{proof}
Our strategy for this theorem is to, given an $x$ sufficiently close to $\xbar$,
\begin{enumerate}
    \item show that the restriction of $\Phi_{f,x}$ along any direction has a strict local minimum, then
    \item argue that among all such local minima, there is one with a smallest value, and
    \item that point is the local minimum of $\Phi_{f,x}$.
\end{enumerate}
For any $x \in \Rn$ and direction $d \in \Rn$ (where without loss of generality we take $\|d\| = 1$), we define the following three functions:
$$a(x,d) \defeq \frac{1}{6}\sum_{i=1}^n d^T d_i \grad_i^3 f(x)d,$$
$$b(x,d) \defeq \frac{1}{2} d^T\Hess f(x)d,$$
$$c(x,d) \defeq \grad f(x)^Td.$$

Then the restriction of $\Phi_{f,x}$ to the line $x + \alpha d$ is the cubic polynomial
$$C_{x,d}(\alpha) \defeq a(x, d)\alpha^3 + b(x, d)\alpha^2 + c(x, d)\alpha + f(x).$$

At $\xbar$, we have for all $d$ $c(\xbar,d) = 0$ for all $d$ and $b(\xbar) \ge \frac{m}{2}$ for some $m$, due to $\xbar$ being a local minimum and $f$ being strongly convex. Now recall that a cubic polynomial $ax^3 + bx^2 + cx$ has a strict local minimum if and only if $b^2 - 3ac > 0$. As $a(x,d), b(x,d),$ and $c(x,d)$ are all continuous in both $x$ and $d$, there exist $\epsilon > 0$ and $R$ such that $b(x, d)^2 - 3a(x, d)c(x, d) > \epsilon$ and $\frac{\sqrt{b(x, d)^2 - 3a(x, d)c(x, d)}}{3a(x,d)} > \epsilon$ whenever $\|x - \xbar\| < R$. As a consequence, $C_{x,d}$ has a local minimum whenever $\|x - \xbar\| < R$.

Now consider any point $x$ satisfying $\|x - \xbar\| < R$. For $d \ne 0$, we define the function
$$\alpha^*(d) \defeq -\frac{c(x,d)}{b(x,d) + \sqrt{b(x,d)^2 - 3a(x,d)c(x,d)}},$$
which is the strict local minimum of $C_{x,d}$. Let
$$\phi^*(d) \defeq C_{x,d}(\alpha^*(d)),$$
which is the value of the local minimum of $C_{x,d}$. By the above conditions, $b(x,d)^2 - 3a(x,d)c(x,d) \ge 0$, and so $\phi^*$ is a continuous function in $d$. As $\|d\| = 1$ is a compact set, we can apply the Weierstrass theorem to guarantee the existence of a global minimum of $\phi^*$, as well as a minimizer $d^*$ of $\phi^*$, which is the direction $d$ along which the local minimum of $C_{x,d}$ has the smallest value.

We now claim that $x + \alpha^*(d^*)d^*$ is a local minimum of $\Phi_{f,x}$. Suppose, for contradiction, that it is not. Then there exists a sequence of points $y_i, i = 1, 2, \ldots$ such that $y_i \to \alpha^*(d^*) d^*$ and $\Phi_{f,x}(x+y_i) < \Phi_{f,x}(x+\alpha^*(d^*)d^*)$. Note that we can use $y_i$ to define sequences $z_i \defeq \frac{y_i}{\|y_i\|}$ and $\alpha_i = \|y_i\|$, so that $y_i = \alpha_iz_i$. Clearly since $y_i \to \alpha^*(d^*)d^*$, it must be that $z_i \to d^*$ and $\alpha_i \to \alpha^*$.

Now let $i$ be sufficiently large so that $|a_i - a^*| < \epsilon/2$. We argue first that 
$$\Phi_{f,x}(x + \alpha_iz_i) > \Phi_{f,x}(x+\alpha^*(z_i)z_i) = \phi^*(z_i).$$
This is because the distance between the inflection point and local minimum of $C_{x,z_i}$ is given by $\frac{\sqrt{b(x,z_i)^2 - 3a(x,z_i)c(x,z_i)}}{3a(x,z_i)}$ (which is by assumption at least $\epsilon$), and all points $\alpha$ satisfying $C_{x,z_i}(\alpha) < C_{x,z_i}(\alpha^*(z_i))$ have to be at least as far away from the local minimum as inflection point is. Now recall by the definition of $d^*$ that $\phi^*(z_i) > \phi^*(d^*)$. We thus have

$$\Phi_{f,x}(x + \alpha_iz_i) > \Phi_{f,x}(x+\alpha^*(z_i)z_i) = \phi^*(z_i) > \phi^*(d^*) = \Phi_{f,x}(x + \alpha^*(d^*)d^*),$$
which is a contradiction.
\end{proof}

We now show that the convergence of this algorithm is commensurate with the order of the derivatives used.

\begin{theorem}\label{Thm: Convergence}
Suppose that a function $f$ is strongly convex with parameter $m$ and has a Lipschitz continuous third derivative with parameter $L$. Suppose further that it has a strict global minimum $x^*$. Then there exist scalars $C$ and $D > 0$ such that if $\|x^0 - x^*\| < D$, the iterates of Algorithm~\ref{Alg: 3ON} satisfy $$\|x^{k+1} - x^*\| \le C\|x^k - x^*\|^3.$$
\end{theorem}
\begin{proof}
Let $x^k$ be the current iterate of Algorithm~\ref{Alg: 3ON}, and suppose that $D$ is such that $D^2 < \frac{m}{L}$ and $x^{k+1}$ is well-defined by Theorem~\ref{Thm: Well-defined}. Further let $\{H_1(x), \ldots, H_n(x)\}$ be the tensor of third derivatives of $f$ at $x$, $Q(x)$ be the Hessian $\Hess f(x)$, and $b(x)$ be the gradient $\grad f(x)$. Then $x^{k+1}$ satisfies

\beq\label{Eq: Gradient zero iterate} \frac{1}{2}\sum_{i=1}^n (x^{k+1}-x^k)_iH_i(x^k)(x^{k+1}-x^k) + Q(x^k)(x^{k+1}-x^k) + b(x^k) = 0, \eeq
and
\beq\label{Eq: Hessian pd iterate} \Hess \Phi_{f,x^k}(x^{k+1}) \succ 0.\eeq
We will further need that $b(x^*) = 0$ and $Q(x^*) \succ 0$, which follow from that $x^*$ is the local minimum of $f$.

We can write

\begin{align}
    0 =& \frac{1}{2}\sum_{i=1}^n (x^{k+1}-x^k)_iH_i(x^k)(x^{k+1}-x^k) + Q(x^k)(x^{k+1}-x^k) + b(x^k) \nonumber\\
    =& \frac{1}{2} \sum_{i=1}^n (x^*-x^k)_iH_i(x^k)(x^*-x^k) + Q(x^k)(x^*-x^k) + b(x^k) \label{Eq: Proof grad line}\\
     &+ \frac{1}{2} \sum_{i=1}^n (x^* - x^k)_iH_i(x^k)(x^{k+1}-x^*) + \frac{1}{2}Q(x^k)(x^{k+1}-x^*)\label{Eq: Proof Hess * line}\\
    &+ \frac{1}{2} \sum_{i=1}^n (x^{k+1} - x^k)_iH_i(x^k)(x^{k+1}-x^*) + \frac{1}{2}Q(x^k)(x^{k+1}-x^*)\label{Eq: Proof Hess k+1 line}
\end{align}

In what follows, we bound each of the quantities (\ref{Eq: Proof grad line}), (\ref{Eq: Proof Hess * line}), and (\ref{Eq: Proof Hess k+1 line}).\\

Observe that (\ref{Eq: Proof grad line}) is $\grad \Phi_{f,x^k}(x^*)$, and also recall that $\grad f(x^*) = 0$. We can thus bound the quantity in (\ref{Eq: Proof grad line}):

\begin{align*}
    0 =& \grad f(x^*)\\
    =& \frac{1}{2} \sum_{i=1}^n (x^*-x^k)_iH_i(x^k)(x^*-x^k) + Q(x^k)(x^*-x^k) + b(x^k)\\
    & +\int_0^1\int_0^1 \sum_{i=1}^n \left(H_i(x^k + \tau\gamma (x^*-x^k))-H_i(x^k)\right)(x^*-x^k)_i(x^*-x^k) \tau d\tau d\gamma.
\end{align*}

The norm of the integral in the second line, and therefore the norm of (\ref{Eq: Proof grad line}), is upper bounded by $\frac{L}{6} \|x^*-x^k\|^3$. 
To see this, note that by (\ref{Lem: Hessian switch}) the summation can be rewritten as
$$d^{\tau,\gamma} \defeq \bmat (x^*-x^k)^TG_1^{\tau,\gamma}(x^*-x^k) \\ \vdots \\ (x^*-x^k)^TG_n^{\tau,\gamma}(x^*-x^k) \emat,$$ where we define $G_i^{\tau,\gamma} \defeq H_i(x^k + \tau\gamma (x^*-x^k))-H_i(x^k)$ for $i \in \otn$. Letting $G^{\tau,\gamma} \defeq \{G_1^{\tau,\gamma}, \ldots, G_n^{\tau,\gamma}\}$, we have by definition,

\begin{align*}
    \|d^{\tau, \gamma}\| =& \max_{\|u\|\le 1} u^Td^{\tau,\gamma}\\
    =& \|x^*-x^k\|^2 \max_{\|u\|\le 1} \sum_{i=1}^n u_i \frac{(x^*-x^k)}{\|x^*-x^k\|}^TG_i^{\tau,\gamma}\frac{(x^*-x^k)}{\|x^*-x^k\|}\\
    \leq& \|x^*-x^k\|^2 \max_{\|u\|, \|v\|, \|w\| \leq 1} \sum_{i=1}^n u_i (v^T G_i^{\tau,\gamma} w) \\
    =& \|x^*-x^k\|^2 \|G^{\tau,\gamma}\|\\
    \le& \tau\gamma L\|x^*-x^k\|^3,
\end{align*}
where the last line follows from the Lipschitz continuity of the third derivative. Hence,

\begin{equation}\label{Eq: Proof norm grad}\|\int_0^1\int_0^1 d^{\tau,\gamma} \tau d\tau d \gamma\| \leq \int_0^1 \int_0^1 \tau^2\gamma L\|x^*-x^k\|^3 d\gamma d\tau \le \frac{L}{6} \|x^*-x^k\|^3.\end{equation}\\

Now observe that (\ref{Eq: Proof Hess * line}) is equal to $\frac{1}{2}\Hess \Phi_{f,x^k}(x^*)(x^{k+1}-x^k)$. We also have
$$\Hess f(x^*) = \Hess \Phi_{f,x^k}(x^*) + \int_0^1\sum_{i=1}^n (x^*-x^k)_i\left(H_i(x^k + \tau (x^*-x^k))-H_i(x^k)\right)d\tau.$$

\noindent Since $\|x^*-x^k\|^2 \le \frac{m}{L}$,
$$\|\int_0^1 \sum_{i=1}^n (x^*-x^k)_i\left(H_i(x^k + \tau (x^*-x^k))-H_i(x^k)\right)d\tau\| \le \frac{m}{2}.$$ 
To see why this is, again let $G_i^\tau \defeq H_i(x^k + \tau (x^*-x^k))-H_i(x^k)$ for $i \in \otn$ and $G^\tau \defeq \{G_1^\tau, \ldots, G_n^\tau\}$, and observe that
\begin{align*}
    \|\sum_{i=1}^n (x^*-x^k)_iG_i^\tau\| =& \max_{\|u\|, \|v\| \leq 1} v^T (\sum_{i=1}^n (x^*-x^k)_iG_i^\tau) u \\
    =& \|x^*-x^k\| \max_{\|u\|, \|v\| \leq 1} \sum_{i=1}^n v^T \frac{(x^*-x^k)_i}{\|x^*-x^k\|} G_i^\tau u \\
    \le& \|x^*-x^k\| \max_{\|u\|, \|v\|, \|w\| \leq 1} \sum_{i=1}^n v^T w_i G_i^\tau u \\
    =& \|x^*-x^k\| \|G^\tau\| \\
    \le& \tau L\|x^*-x^k\|^2,
\end{align*} where the last line follows from the Lipschitz continuity of the third derivative. Hence, $$\|\int_0^1 \sum_{i=1}^n (x^*-x^k)_iG_i^\tau d\tau \| \leq \int_0^1 \tau L\|x^*-x^k\|^2 d\tau \le \frac{L}{2} \|x^*-x^k\|^2 \le \frac{m}{2}.$$

We then get

\begin{align*} \|\sum_{i=1}^n (x^* - x^k)_iH_i(x^k) + Q(x^k)\| =& \|\Hess f(x^*) - \int_0^1 \sum_{i=1}^n (x^*-x^k)_i\left(H_i(x^k + \tau (x^*-x^k))-H_i(x^k)\right)d\tau\|\\
\ge & \|\Hess f(x^*)\| - \|\int_0^1 \sum_{i=1}^n (x^*-x^k)_i\left(H_i(x^k + \tau (x^*-x^k))-H_i(x^k)\right)d\tau\|\\
\ge& m - \frac{m}{2} = \frac{m}{2}.
\end{align*}
Note in particular that $\sum_{i=1}^n (x^* - x^k)_iH_i(x^k) + Q(x^k)$ is a positive definite matrix whose smallest eigenvalue is at least $\frac{m}{2}$. 

Finally, observe that (\ref{Eq: Proof Hess k+1 line}) is $\frac{1}{2}\Phi_{f,x^k}(x^{k+1})(x^{k+1}-x^*)$, where $\Phi_{f,x^k}(x^{k+1})$ is a positive definite matrix. Putting everything together, we have have $-(\ref{Eq: Proof grad line}) = (\ref{Eq: Proof Hess * line}) + (\ref{Eq: Proof Hess k+1 line})$, which yields
\begin{align*}
x^{k+1}-x^* =& -\left(\frac{1}{2}(\sum_{i=1}^n (x^*-x^k)_iH_i(x^k) + Q(x^k)) + \frac{1}{2}(\sum_{i=1}^n (x^{k+1}-x^k)_iH_i(x^k) + Q(x^k))\right)^{-1}\\
& \times \left(\frac{1}{2} \sum_{i=1}^n (x^*-x^k)_iH_i(x^k)(x^*-x^k) + Q(x^k)(x^*-x^k) + b(x^k)\right).
\end{align*}
The first factor of the right hand side is the inverse of a positive definite matrix whose minimum eigenvalue is at least $\frac{m}{4}$, and thus is a matrix whose norm is at most $\frac{4}{m}$. The second factor is a vector with norm at most $\frac{L}{2}\|x^*-x^k\|^3$ (from (\ref{Eq: Proof norm grad})). Hence we find
$$\|x^{k+1} - x^*\| \le \frac{L}{6m}\|x^k-x^*\|^3$$
as desired.
\end{proof}

We end this section by establishing that whenever the third-order Taylor approximation has a strict local minimum, (\ref{Eq: complete cubic SDP}) is strictly feasible. By strict feasibility, we mean that there exists $(X,x,y)$ such that the equality constraints are satisfied and the semidefinite constraint holds with strict positive definiteness. While this condition is not necessary for (\ref{Eq: complete cubic SDP}) to have a solution, it is nonetheless a relevant computational consideration for interior point methods.

\begin{theorem}\label{Thm: SLM strictly feasible}
    If the cubic polynomial associated with (\ref{Eq: complete cubic SDP}) has a strict local minimum, then (\ref{Eq: complete cubic SDP}) is strictly feasible.
\end{theorem}

\begin{proof}
Let $p$ be the cubic polynomial associated with (\ref{Eq: complete cubic SDP}) and let $\xbar$ be its strict local minimum. Note that since $y$ is otherwise unconstrained, for any matrix $X$ and vector $x$ satisfying $\Hess p(x) \succ 0$, there exists $y$ large enough such that $$\bmat \sum_{i=1}^n x_iH_i + Q & \sum_{i=1}^n \Tr(H_iX)e_i+Qx \\ (\sum_{i=1}^n \Tr(H_iX)e_i+Qx)^T & y\emat \succ 0.$$ This is because the matrix in the top-left-hand corner is in fact $\Hess p(x)$. Therefore it suffices to construct $(X,x)$ satisfying the equality constraints such that $\Hess p(x) \succ 0$ and $X - xx^T \succ 0$.

Let $D \in \R^{n \times n}$ be any positive definite matrix, and let $d$ be the vector in $\Rn$ where the $i$-th entry is $\Tr(H_iD)$. Finally, define $$v \defeq \frac{1}{2}(\sum_{i=1}^n \xbar_iH_i + Q)^{-1}d.$$

Consider scaling $D$ by a positive scalar $\alpha$, which also scales $v$ by $\alpha$. For $\alpha$ sufficiently small, we have
$$\Hess p(\xbar - \alpha v) \succ 0,$$
since $\Hess p(\xbar - \alpha v) \to \Hess p(\xbar)$ as $\alpha \to 0$. Similarly, for (a potentially smaller) sufficiently small $\alpha > 0$, we have
$$\alpha D - \alpha^2 vv^T = \alpha (D - \alpha vv^T) \succ 0.$$
For the smaller of these two scalings, let $\ybar \defeq \xbar - \alpha v$. Without loss of generality, suppose that $D$ was such that $\alpha = 1$, that is, $D$ was chosen to be small enough that no scaling was necessary.

We now set
$$X \defeq \ybar\ybar^T + D - vv^T, x = \ybar.$$

By construction, $\Hess p(\ybar) \succ 0$. We know that the third constraint of (\ref{Eq: complete cubic SDP}) holds with positive definiteness from the Schur complement condition and that $X - xx^T = D - vv^T \succ 0$. It only remains to verify the equality constraints
$$\frac{1}{2}\Tr(H_iY) + e_i^TQy + b_i = 0, \forall i.$$
 This is proven by the following chain of equalities:
\begin{align*}
    \frac{1}{2}\Tr(H_iX) + e_i^TQx + b_i &= \frac{1}{2} \Tr(H_i\ybar\ybar^T) + \frac{1}{2}\Tr(H_iD) - \frac{1}{2}\Tr(H_ivv^T) + e_i^TQx + b_i\\
    &= \frac{1}{2}(\xbar - v)^TH_i(\xbar - v) + \frac{1}{2}\Tr(H_iD) - \frac{1}{2}v^TH_iv + e_i^TQ(\xbar - v) + b_i\\
    &=\frac{1}{2}\xbar^TH_i\xbar + e_i^TQ\xbar + b_i - \xbar^TH_iv - e_i^TQv + \frac{1}{2}\Tr(H_iD)\\
    &= 0 -\left(\sum_{i=1}^n (\xbar_iH_i + Q)v\right)_i + \frac{1}{2}\Tr(H_iD)\\
    &= -\frac{1}{2}d_i + \frac{1}{2}\Tr(H_iD)\\
    &= 0.
\end{align*}

Between the third and fourth lines, we used that $\xbar$ was a stationary point, and then applied the definition of $v$ between the fourth and fifth lines.

\end{proof}

\section{Adaptations for Infeasible Iterations}\label{Sec: Adaptations}

Though we have established that our algorithm is well-defined when close to the global minimum, we should nonetheless not expect this to be the case. Rather, for most initializations it is likely that the cubic approximation has no local minimum. Hence, in this section we present two revised iteration strategies in the event the problem is not feasible.

\subsection{A ``Convex'' Problem}\label{SSec: Convex Problem}

We first consider the following problem, which is more fully discussed in Section 4.2 of \cite{ahmadi2022complexity}:
\begin{equation}\label{Eq: minimize cubic}
    \begin{aligned}
	& \underset{x^{k+1} \in \Rn}{\inf} 
	& & \Phi_{f,x^k}(x^{k+1})\\
	& \text{subject to}
	&& \Hess \Phi_{f,x^k}(x^{k+1}) \succeq 0,\\
    \end{aligned}
\end{equation}
which has the interpretation of minimizing $\Phi_{f,x^k}(x^{k+1})$ over the region where $\Phi_{f,x^k}(x^{k+1})$ is convex, which is itself a convex set. In the event $\Phi_{f,x^k}(x^{k+1})$ has a local minimum, the solutions of (\ref{Eq: minimize cubic}) and (\ref{Eq: unrelaxed problems}) coincide. Though (\ref{Eq: unrelaxed problems}) does not have a solution when no local minimum exists, (\ref{Eq: minimize cubic}) is nonetheless bounded below when $\Hess p(x^{k}) \succ 0$.

\begin{proposition} If $\Hess f(x^{k}) \succ 0$, then (\ref{Eq: minimize cubic}) is bounded below and attains its solution.\end{proposition}

\begin{proof}

Note that by assumption the feasible set is non-empty. We show that the minimizer of (\ref{Eq: minimize cubic}) must lie within a bounded region. We define similarly as in the proof of Theorem~\ref{Thm: Well-defined}:
$$a(d) \defeq \frac{1}{6}\sum_{i=1}^n d^T d_i \grad_i^3 f(x^k)d,$$
$$b(d) \defeq \frac{1}{2} d^T\Hess f(x^k)d,$$
$$c(d) \defeq \grad f(x^k)^Td.$$
$$C_{d}(\alpha) \defeq a(d)\alpha^3 + b(d)\alpha^2 + c(d)\alpha + f(x^k).$$

Along any direction $d$, the function $C_d$ has either an inflection point or a local minimum (possibly both). The inflection point if it exists, occurs at $\frac{-b(d)}{3a(d)}$, and the local minimum, if it exists, occurs at $-\frac{c(d)}{b(d) + \sqrt{b(d)^2 - 3 a(d)c(d)}}$. Note that at least one of the quantities is well-defined for any given $d$.

Consider now the quantity

\begin{equation}
    \begin{aligned}
	R \defeq& \underset{d \in \Rn, \|d\|=1}{\sup} 
	& & \min\left(\left|\frac{b(d)}{3a(d)}\right|, \left|\frac{c(d)}{b(d) + \sqrt{b(d)^2 - 3 a(d)c(d)}}\right|\right),\\
    \end{aligned}
\end{equation}

We claim that if $R$ exists, then the minimum of (\ref{Eq: minimize cubic}) must satisfy $\|x^{k+1} - x^k\| \le R$. To see why, first recall that for $\Hess \Phi_{f,x^k}(x) \succeq 0$ to hold at some $x$, the second derivative must be nonnegative along all restrictions to a line. Hence the feasible set of (\ref{Eq: minimize cubic}) cannot extend beyond the inflection point along any direction. Further recall that the region $\{x | \Hess \Phi_{f,x^k}(x) \succeq 0\}$ is convex and contains $x^k$. In particular, this implies that no point that extends beyond the local minimum in any direction can be the the solution to (\ref{Eq: minimize cubic}), as the local minimum along that direction is feasible.

It remains to show that $R$ exists. To see why, first observe that whenever $b(d)^2 - 3 a(d)c(d) \ge 0$,
$$\left|\frac{c(d)}{b(d) + \sqrt{b(d)^2 - 3 a(d)c(d)}}\right| \le \frac{|c(d)|}{b(d)} \le \frac{\|\grad f(x^k)\|}{\lambda},$$
where $\lambda$ is the smallest eigenvalue of $\Hess f(x^k)$ (which is positive by assumption).

If instead $b(d)^2 - 3 a(d)c(d) < 0$, then we have
$$a(d) \ge \frac{\lambda^2}{\|\grad f(x^k)\|} \Rightarrow \left|\frac{b(d)}{3a(d)}\right| \le \frac{\|\grad f(x^k)\|}{3\lambda}.$$
Thus $R$ is finite. Since the objective of (\ref{Eq: minimize cubic}) is continuous and the feasible set is closed and contained in a compact ball of radius $R$, the minimum is attained.
\end{proof}

One can show that a valid lifting relaxation of (\ref{Eq: minimize cubic}) is given by
\begin{equation}\label{Eq: no gradient cubic SDP}
    \begin{aligned}
	& \underset{X \in \Snn, x \in \Rn, y \in \R}{\inf}
	& & \frac{1}{3}\Tr(QX) + b^Tx + \frac{y}{6}\\
	& \text{subject to}
	&& \bmat \sum_{i=1}^n x_iH_i + Q & \sum_{i=1}^n \Tr(H_iX)e_i+Qx \\ (\sum_{i=1}^n \Tr(H_iX)e_i+Qx)^T & y\emat \succeq 0,\\
	&&& \bmat X & x \\ x^T & 1 \emat \succeq 0,
    \end{aligned}
\end{equation}
which is the same as (\ref{Eq: cubic SDP}) except without the lifted gradient constraints. In the event that $\Hess p(x^k) \succeq 0$ (which, for example, is always the case if $f(x)$ is convex), this SDP is guaranteed to be feasible, as the (lifted) origin itself is a feasible solution. Using this formulation, after extracting the the vector $x$ from the solution, $x^k + x$ serves as the next iterate.

\subsection{Levenberg-Marquardt Regularization for the Third-Order Method}\label{SSec: LM Regularization}

In the classical Newton method for nonconvex optimization, Levenberg–Marquardt (LM) regularization adds a multiple of the identity matrix to the Hessian to ensure positive definiteness. This guarantees that the resulting quadratic approximation has a well-defined global minimum and stabilizes the iterates when far from optimality.

We adopt a similar idea for the third-order method by adding a sufficiently large multiple of the identity matrix to the Hessian at the current iterate. Equivalently, this corresponds to adding a regularization term of the form $\alpha \|x - \xbar\|_2^2$ to the objective. We emphasize that this differs from existing convex-regularized higher-order methods: here, the regularization term has \emph{lower} degree than the highest-order Taylor term. Our goal is not to enforce global convexity, but rather to guarantee the existence of a local minimum. Importantly, this preserves the cubic structure of the model and thus its tractability via semidefinite programming.

It is not entirely trivial to compute the ideal parameter $\alpha$. One cannot for example solve an SDP feasibility problem to determine the smallest value ensuring the existence of a local minimum, as this would induce a quadratic semidefinite constraint. One could instead use an increasing sequence of values until feasibility is achieved, but this requires solving additional SDPs. Instead, we propose an explicit choice of $\alpha$ computable directly from the derivatives.

For an iterate $x$, let $\lambda$ be the minimum eigenvalue of $\Hess f(x)$ and $$g \defeq \bmat |\grad_1 f(x)| \\ \vdots \\ |\grad_n f(x)| \emat ,h \defeq \bmat \|\grad_1^3 f(x)\| \\ \vdots \\ \|\grad_n^3 f(x)\| \emat.$$
We set
\beq \label{Eq: alpha LM} \alpha_{LM} \defeq \sqrt{\frac{3}{2}(\|g\|\|h\| + g^Th)} - \min \{0,\lambda\}.\eeq

\begin{proposition} For a three-times differentiable function $f$ and any point $\xbar \in \Rn$, the function $\Phi_{f,\xbar}+\alpha_{LM}\|x-\xbar\|^2$ has a local minimum.\end{proposition}

\begin{proof}
Recall from the proof of Theorem~\ref{Thm: Well-defined} that if the function \beq  \delta(x,d) \defeq b(x,d)^2 - 3a(x,d)c(x,d) \eeq as defined in the proof is positive for all $d \ne 0$, then the Taylor expansion of $f$ around $x$ has a local minimum. If for some $\alpha \in \R$, $\Hess f(x)$ is replaced by $\Hess f(x) + \alpha I$, then this expression becomes
\beq \label{Eq: nonneg condition} (b(x,d) + \alpha d^Td)^2 - 3a(x,d)c(x,d).\eeq
Since this function is still homogeneous with respect to $d$, we can consider whether this function is positive on $\|d\| = 1$. If we can find $\beta$ such that $a(x,d)c(x,d) \le \beta, \forall \|d\| = 1$, then we can set
\beq \label{Eq: norm alpha} \alpha = \sqrt{3\beta} - \min\{0, \lambda\},\eeq
and satisfy (\ref{Eq: nonneg condition}). We point out that one natural choice of $\beta$ is $\|\grad^3 p(x)\| \|\grad p(x)\|$, but it is NP-hard to compute the norm of the tensor $\grad^3 p(x)$~\cite{hillar2013most}. Instead, for the rest of the proof we establish that \beq\label{Eq: 3ac bound} \beta = \frac{\|g\|\|h\| + g^Th}{2}\eeq
is a valid choice for $\beta$. Recall that
$$a(x,d)c(x,d) = \left(\grad f(x)^Td\right)\left(\sum_{i=1}^n d_i (d^T\grad_i^3 f(x)d)\right).$$

One can see that

$$\left|a(x,d)c(x,d)\right| = \left|(\grad f(x)^Td)(\sum_{i=1}^n d_i (d^T\grad_i^3 f(x)d))\right| \le \left(\sum_{i=1}^n |d_i|g_i\right)\left(\sum_{i=1}^n |d_i|h_i\right).$$
For the purpose of analyzing the final expression on the right hand side, without loss of generality, we consider $d$ in the nonnegative orthant. For such $d$, since $g$ and $h$ are nonnegative vectors, the final expression then reduces to $(d^Tg)(h^Td).$

Now observe that
$$(d^Tg)(h^Td) = d^T \left(\frac{gh^T + hg^T}{2}\right)d,$$ where $gh^T + hg^T$ is a rank-2 matrix with the following eigenvalue-eigenvector pairs:
$$\left(g^Th + \|g\|\|h\|, \frac{g}{\|g\|} + \frac{h}{\|h\|}\right)$$

$$\left(g^Th - \|g\|\|h\|, \frac{g}{\|g\|} - \frac{h}{\|h\|}\right)$$
To verify this, one can see that

\baeq\left(gh^T + hg^T\right)\left(\frac{g}{\|g\|} + \frac{h}{\|h\|}\right) =& \frac{g}{\|g\|}h^Tg + \frac{g}{\|h\|}h^Th + \frac{h}{\|g\|}g^Tg + \frac{h}{\|h\|}g^Th\\
=&\left(\frac{g}{\|g\|} + \frac{h}{\|h\|}\right)g^Th + \frac{\|g\|}{\|g\|}\|h\|g + \frac{\|h\|}{\|h\|}\|g\|h\\
=& (g^Th + \|g\|\|h\|)\left(\frac{g}{\|h\|} + \frac{h}{\|h\|}\right) &\eaeq
The other pair can be shown analogously. Note that because $g$ and $h$ have nonnegative entries by definition, the larger of the two eigenvalues is $g^Th + \|g\|\|h\|$. Substituting this bound into (\ref{Eq: norm alpha}) yields $\alpha_{LM}$, completing the proof.
\end{proof}

It should be noted that $\alpha_{LM}$ will in general overestimate the regularization necessary to ensure the feasibility of (\ref{Eq: cubic SDP}). This may be undesirable as larger regularization terms lead to less fidelity in the approximation and smaller steps and convergence. Our recommended approach is to combine ideas from both adaptations. One can choose the minimum eigenvalue of $\Hess f(x)$ as the regularization constant (plus a small constant), and apply the SDP (\ref{Eq: no gradient cubic SDP}) instead. This yields a regularization strategy that aligns itself analogously with that of Levenberg-Marquardt for second-order models.

\section{Simulations}\label{Sec: Simulations}

In this section, we provide simulation studies. We omit studies benchmarking our method against others; instead we refer the reader to \cite{cai2026globally}. The studies we perform in this work will be focused comparing the infeasible iteration strategies in Section~\ref{Sec: Adaptations} (which were not available to the authors of \cite{cai2026globally}) and the qualitative comparison of Newton fractals.

\begin{figure}[H]
    
    \centering
        \begin{subfigure}[b]{0.48\linewidth}
        \includegraphics[width=\linewidth]{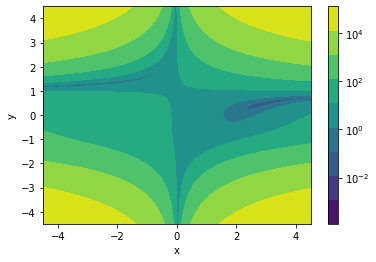}
        \caption{Contour plot of the Beale function.}
    \end{subfigure}
    \begin{subfigure}[b]{0.48\linewidth}
        \includegraphics[width=\linewidth]{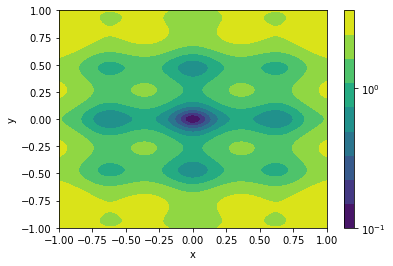}
        \caption{Contour plot of the Bohachevsky function.}
    \end{subfigure}
    \begin{subfigure}[b]{0.48\linewidth}
        \includegraphics[width=\linewidth]{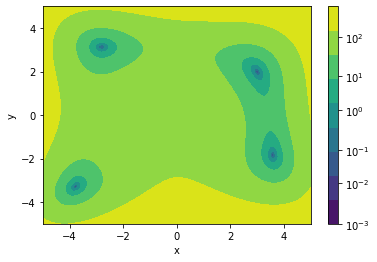}
        \caption{Contour plot of the Himmelblau function.}
    \end{subfigure}
    \begin{subfigure}[b]{0.48\linewidth}
        \includegraphics[width=\linewidth]{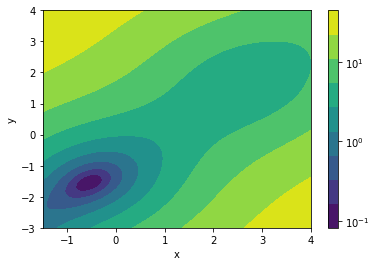}
        \caption{Contour plot of the McCormick function.}
    \end{subfigure}
    \caption{Contour plots of our selected functions.}
    \label{Fig: Contour plots}
\end{figure}

For our experiments, we chose the Bohachevsky, McCormick, Beale, Himmelblau, and Rosenbrock functions, which are some benchmark functions which appear in the literature~\cite{Adorio2005MVFM}. All their function definitions are in Table~\ref{Tab: Function defs}, and their contour plots are in Figure~\ref{Fig: Contour plots}.

\begin{table}[H]
\centering
\begin{tabular}{|l| l  | c |} \hline
Function & Definition & Global minimum \\\hline
Beale & $(1.5-x+xy)^2 + (2.25 - x + xy^2) + (2.625 - x + xy^3)^2$ & $(3,.5)$\\\hline
Bohachevsky & $x^2 + 2y^2 - .3\cos(3\pi x) - .4\cos (4\pi y) + .7$ & $(0,0)$\\\hline
Himmelblau & $(x^2 + y - 11)^2 + (x + y^2 - 7)^2$ & (3,2)\\\hline
McCormick & $\sin (x+y) + (x-y)^2 - 1.5x + 2.5y + 1$ & None \\\hline
\end{tabular}
\caption{Definitions of our selected functions.}
\label{Tab: Function defs}
\end{table}

The Beale function has a cross-shaped valley, with peaks in the corners of the $x,y$-plane and a global minimum at $(3,.5)$. The Bohachevsky function is characterized with an overall bowl shape which has a global minimum at $(0,0)$. However, this function has numerous local minima at intervals of approximately $.61$ in the $x$ direction and $.46$ in the $y$ direction, making it very easy for optimization algorithms to get stuck in a local minimum. Finally, the Himmelblau function is a bowl-shaped function, with four local minima: $(3,2), (-2.805118,3.131312), (-3.779310,-3.283186),$ and $(3.584428,-1.848126)$. Out of these, $(3,2)$ is the global minimum. The McCormick is a plate-shaped function, with no global minimum but many local minima at intervals of approximately $(\pi, \pi)$ of each other. These local minima have progressively smaller objective values as $x$ and $y$ approach $-\infty$. The usual search interval of this function is $-1.5 \le x \le 4, -3 \le y \le 4$, over which the minimum is $(-.54719, -1.54719)$.

In our experiments, we will compare the performance with two objective functions and two regularization strategies. This yields four different strategies which we will compare using Dolan-More performance profiles.

For the objective function, we will compare the objective functions discussed in (\ref{Eq: complete cubic SDP}) and (\ref{Eq: cubic SDP}. We will refer to these as the ``gap'' and ``cubic'' objectives respectively.

For the regularization strategy, we will compare using a heuristic Levenberg-Marquardt approach that uses $\alpha_{LM}$ as discussed in Section~\ref{SSec: LM Regularization} versus one that enacts (\ref{Eq: no gradient cubic SDP}), adding only a regularization parameter that makes the Hessian positive semidefinite if needed. We will refer to these as the Levenberg-Marquardt and heuristic strategies respectively.

For each strategy, we initialize the algorithm in its given search region on a grid of 101 $\times$ 101 starting points and limit to at most 100 iterations. This yields a total of $101 \times 101 \times 4 = 40,804$ initializations, and from each we track the number of iterations required.

\subsection{Algorithm Specification Comparisons}\label{SSec: Comparisons}

Below is the Dolan–Moré~\cite{dolan2002benchmarking} landscape for the four strategies. For a given initialization of test problems, we measure the number of iterations needed to converge. For each problem, the best-performing method is assigned a baseline value, and all other methods are compared relative to this baseline via a performance ratio. The Dolan–Moré profile then plots, for each method, the fraction of problems for which its performance ratio is within a factor $\tau$ of the best. The value at $\tau = 1$ represents the proportion of problems on which a method is the best, while larger values of $\tau$ indicate how often a method performs within a given multiple of the best. Asymptotically, the curve represents the proportion of problems for which the method converged relative to the others.

A few patterns emerge in the Dolan–Moré landspace. Notably, the third-order method converges in fewer iterations and with greater reliability than the second-order method does. For lower values of $\tau$, the specifications using the LM regularization converged in fewer iterations on average. For specifications using the same regularization strategy, those that used the cubic objective function converged in fewer iterations. For higher values of $\tau$, the main separating factor is the choice of objective function, where the cubic objective demonstrates better reliability than the gap objective. However, the difference asymptotically is minimal.

\begin{figure}[H]
    \centering
    \includegraphics[]{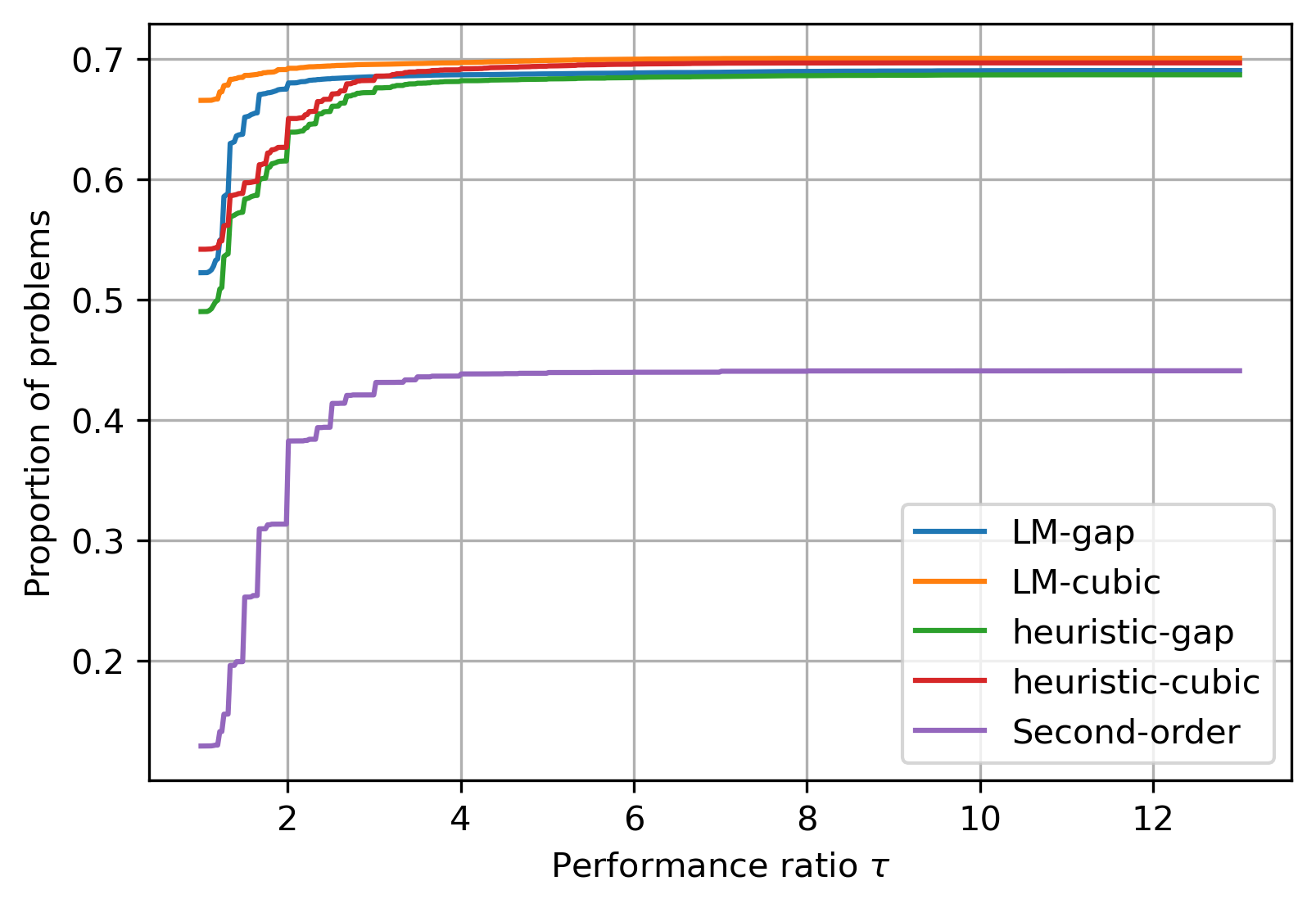}
    \caption{Dolan–Moré of four specifications for the third-order Newton algorithm. Specifications are compared on a total of 40,804 initializations. The curves represent the proportion of instances for which a given specification performed within a factor $\tau$ of the best specification.}
    \label{fig:placeholder}
\end{figure}

\subsection{Newton Fractals}\label{SSec: Newton Fractals}
In this subsection, we compare Newton fractals of the third order method to the second order method (with Marquardt-Levenberg regularization only). A Newton fractal is a visual way to analyze the sensitivity of an algorithm to the initial starting point. Each pixel represents an initial starting point for one of the Newton algorithms. If two pixels have the same color, then the Newton algorithm converges to the same point when starting from that point. Generally speaking, it is desirable when regions whose initial points converge to the same limit are contiguous, and do not display ``fractal'' behavior. For the most part, the Newton fractals from the third-order method indicate less fractal behavior than those arising from the second-order method, and have fewer isolated regions that seem to converge to a particular global minimum by coincidence.

Since most of the third-order fractals are similar, we refer the reader to Appendix~\ref{Sec: Complete Newton Fractals} for the full list of Newton Fractals; in this section we only compare and contrast the third-order Newton fractals for the cubic objective and the regularization.

\begin{figure}[H]
    \centering

    \begin{tabular}{c c c}
        & \textbf{Second-order} & \textbf{Third-order} \\

        \raisebox{7em}{\rotatebox[origin=c]{90}{\textbf{Without Regularization}}} &
        \includegraphics[width=0.45\textwidth]{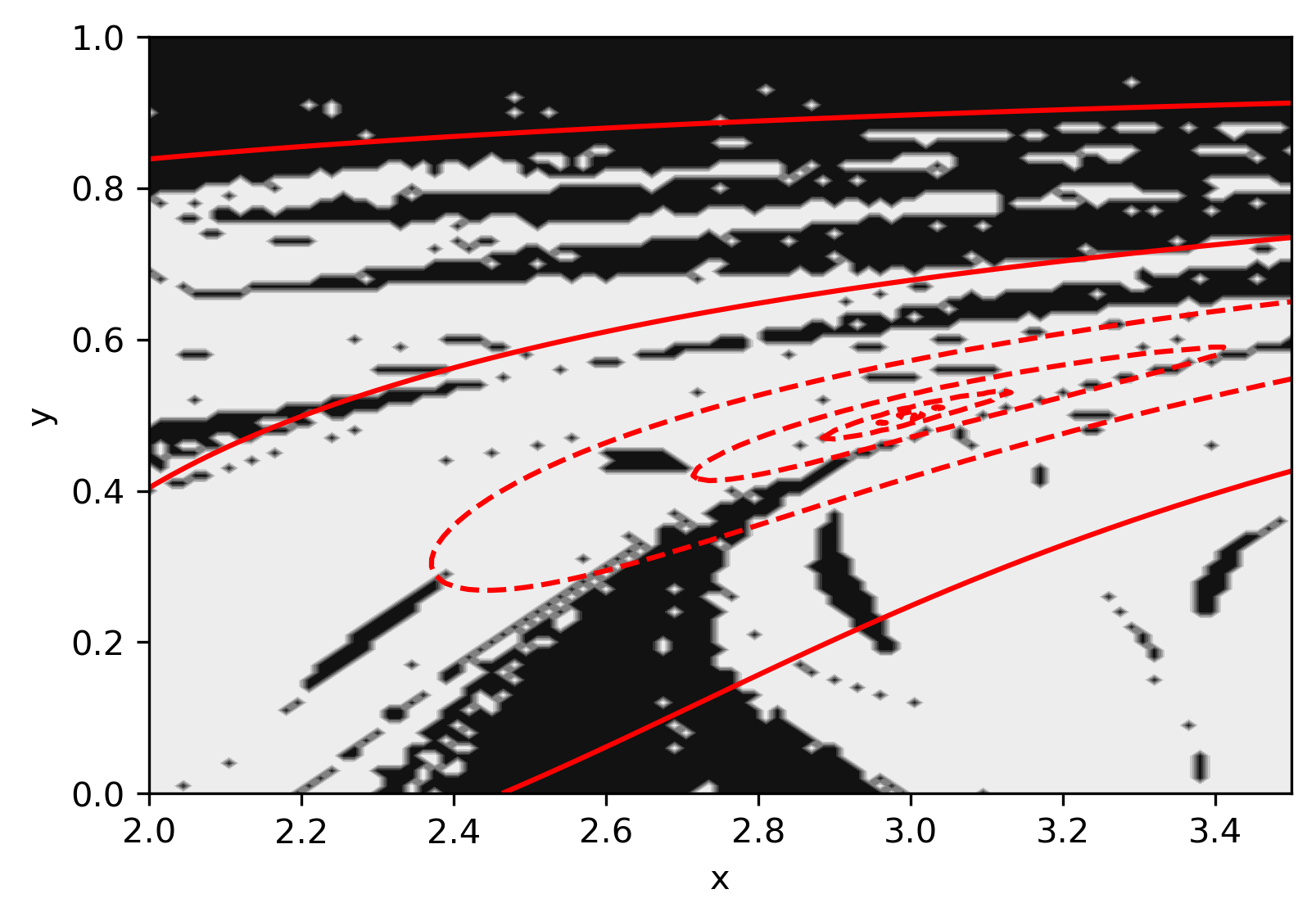} &
        \includegraphics[width=0.45\textwidth]{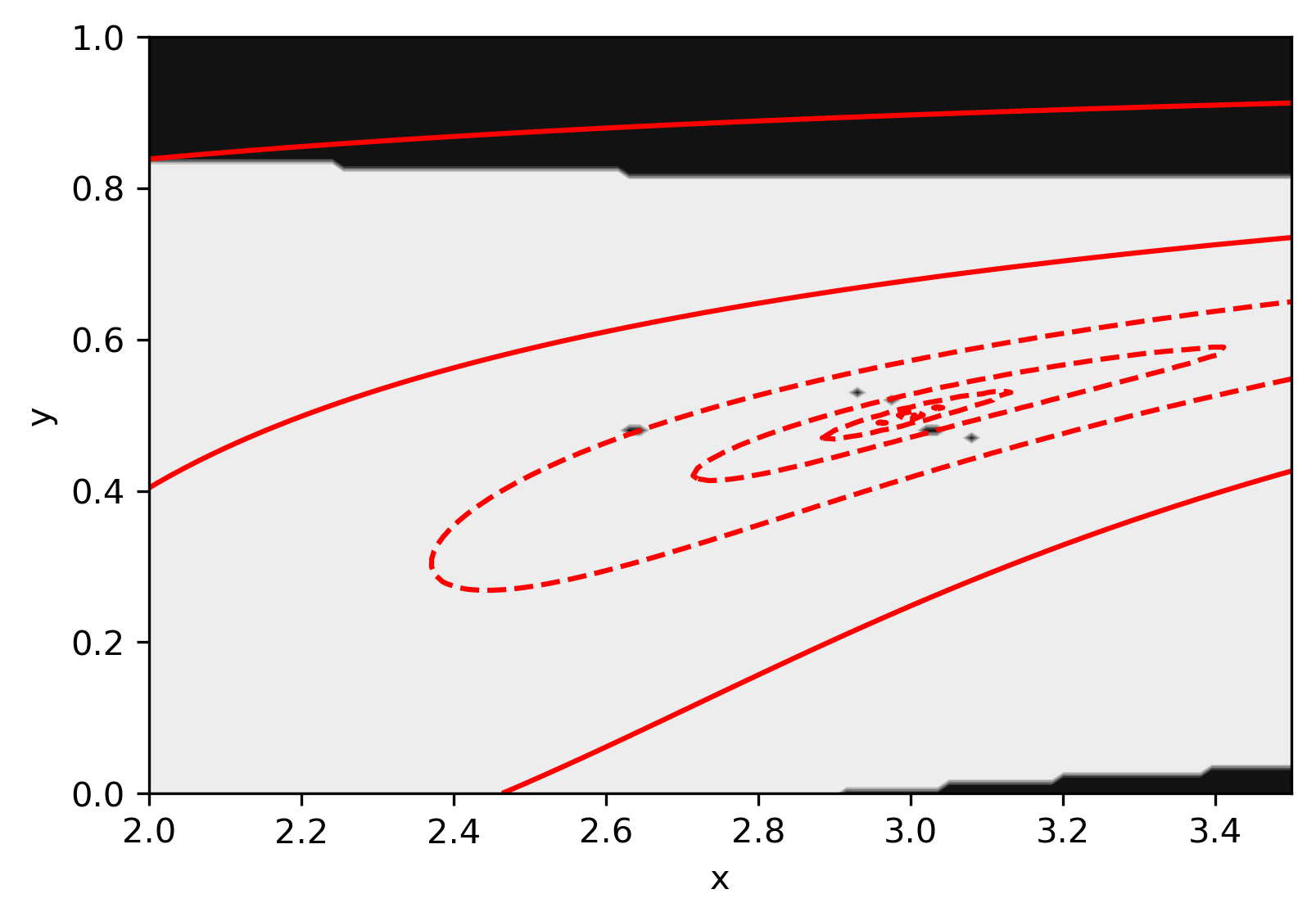} \\

        \raisebox{7em}{\rotatebox[origin=c]{90}{\textbf{With Regularization}}} &
        \includegraphics[width=0.45\textwidth]{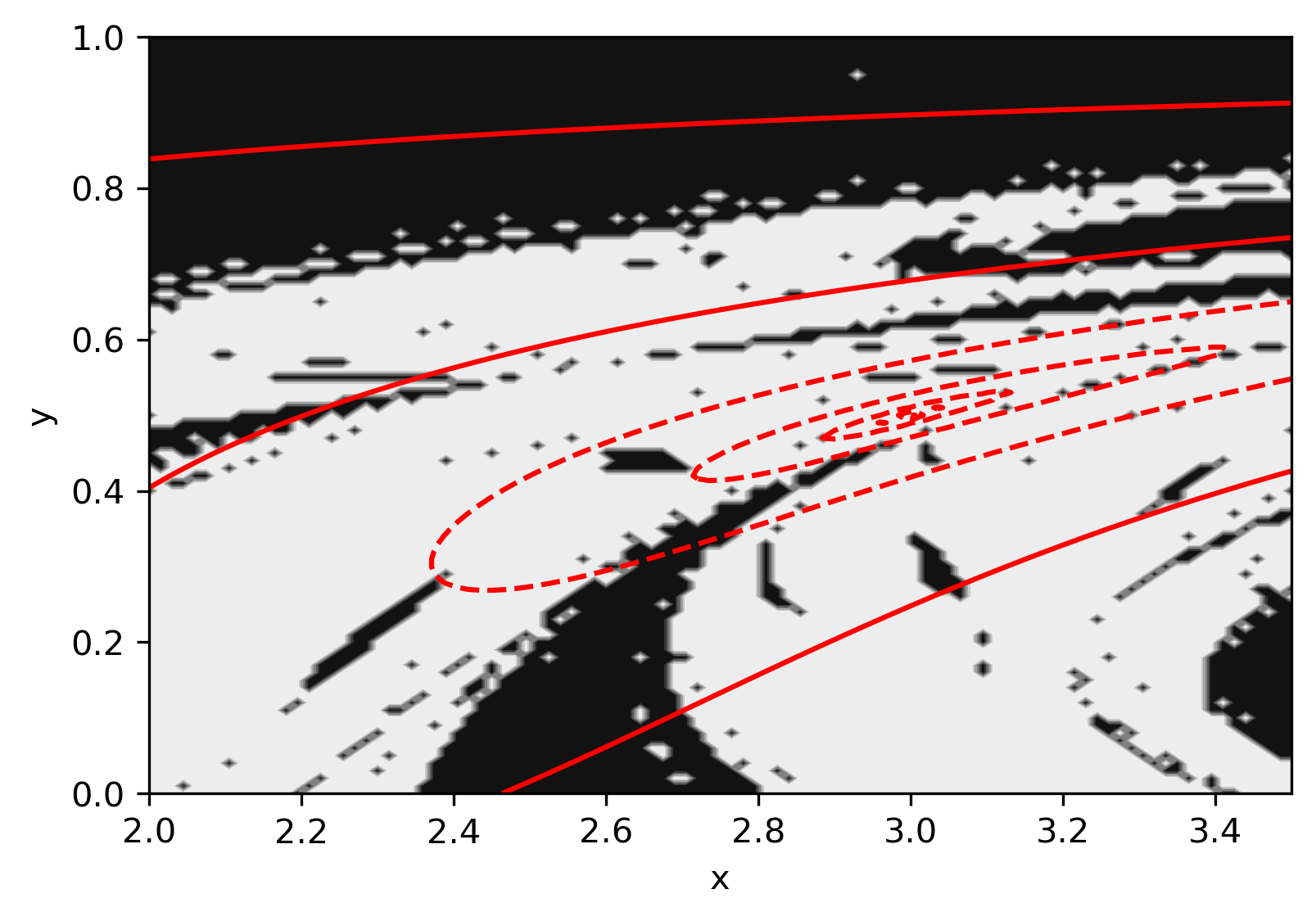} &
        \includegraphics[width=0.45\textwidth]{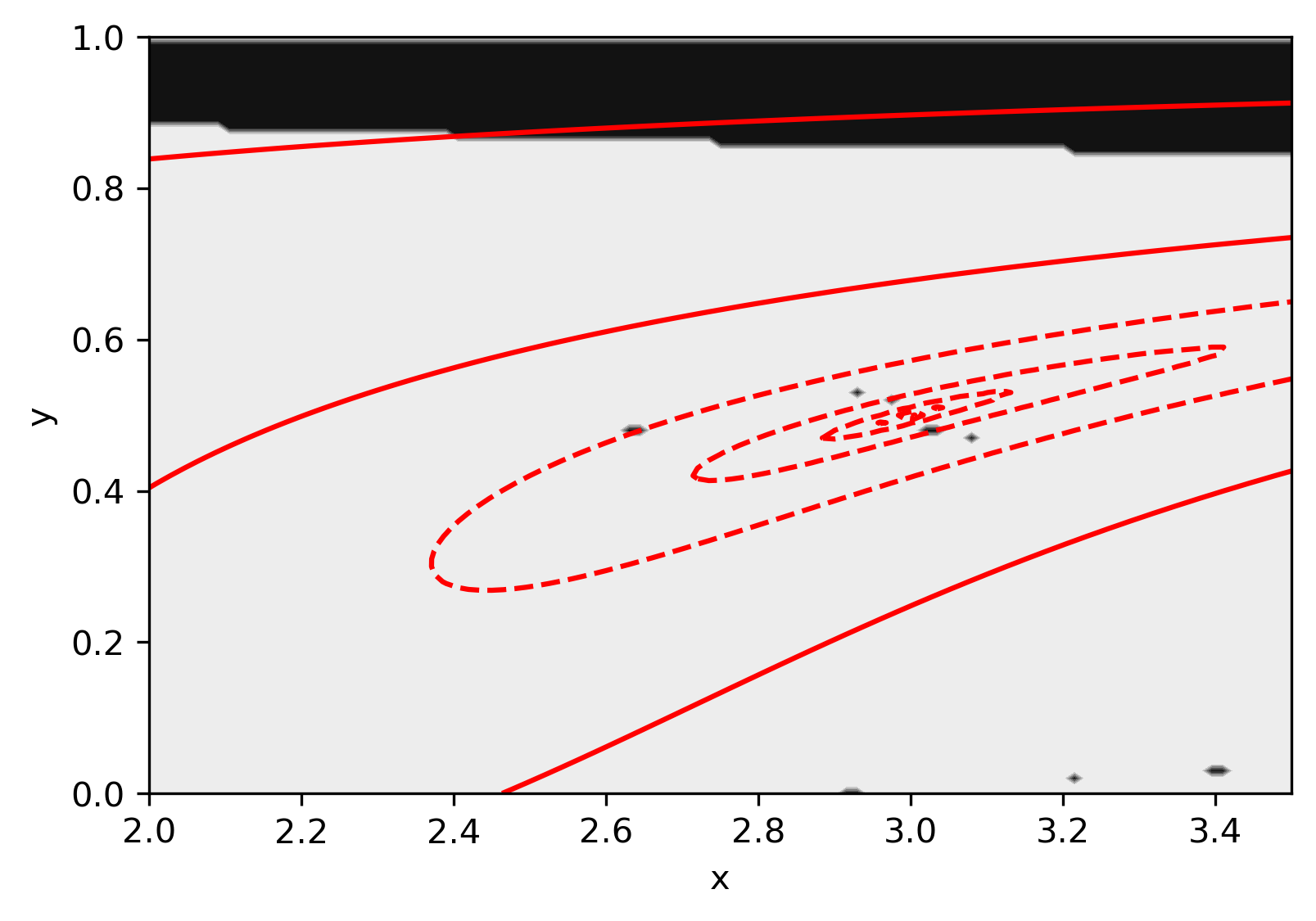} \\
    \end{tabular}

    \caption{
    Comparison of basins of attraction for second- and third-order methods for the Beale function, with and without regularization. Black regions indicate initializations that did not converge within 100 iterations. Otherwise, differently colored regions indicate initial starting points which converged to the same point.
    }
    \label{Fig: Beale Fractals}
\end{figure}

The fractal for the second-order method demonstrates very fractal behavior, with many iterates converging to a saddle point located at (0, 1). By comparison, the third-order fractal displays very stable behavior, with the basin of attraction being a contiguous region around the global minimum. Of the fractals considered, this is is the only function whose appearance is sensitive to the iteration budget. The initializations upper region that are colored black are convergent, but do not do so within 100 iterations. This is more apparent in Figure~\ref{Fig: Beale Fractals 3} which compares the fractals for different specifications of the third-order method.

\begin{figure}[H]
    \centering

    \begin{tabular}{c c c}
        & \textbf{Second-order} & \textbf{Third-order} \\

        \raisebox{7em}{\rotatebox[origin=c]{90}{\textbf{Without Regularization}}} &
        \includegraphics[width=0.45\textwidth]{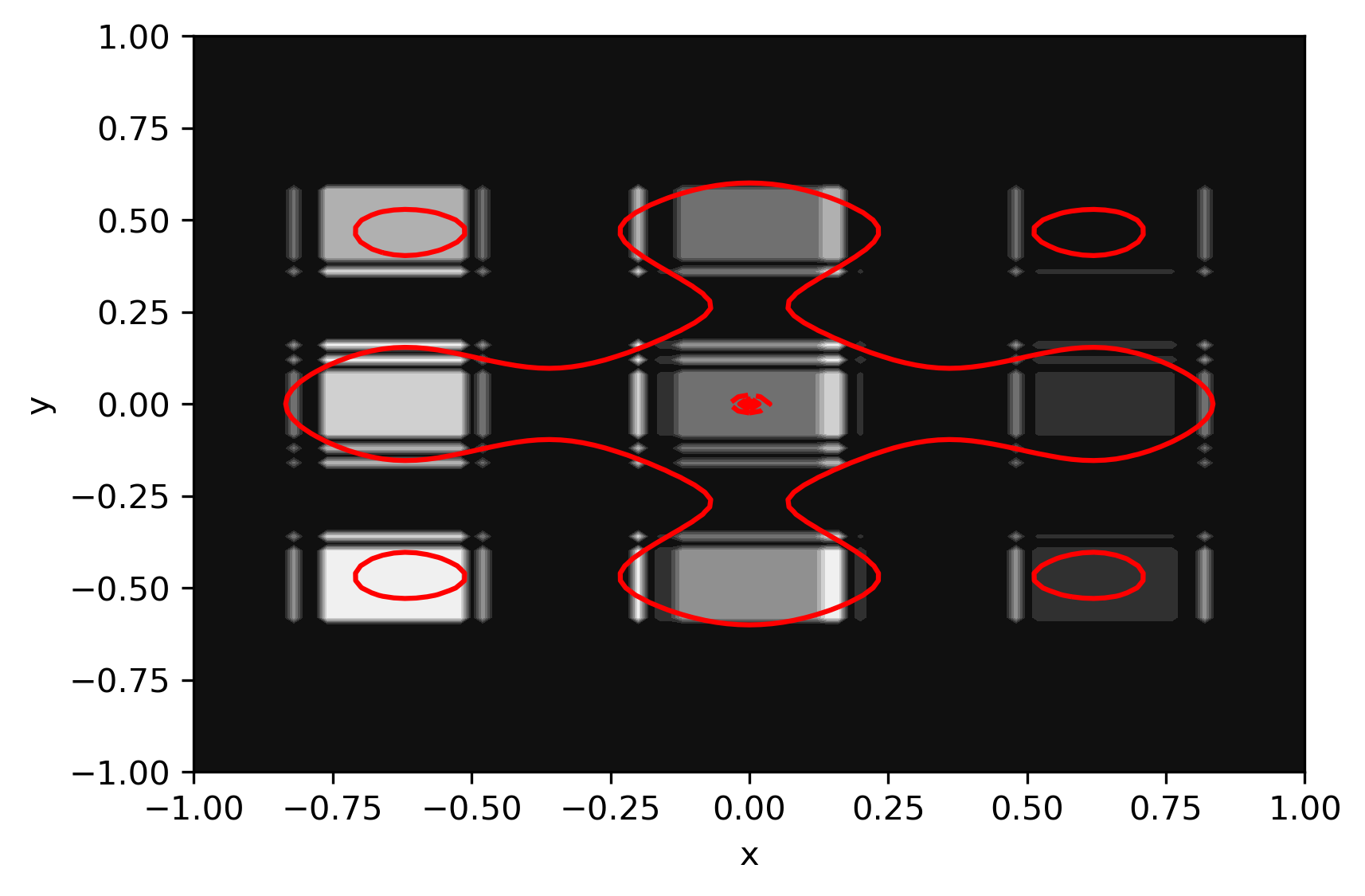} &
        \includegraphics[width=0.45\textwidth]{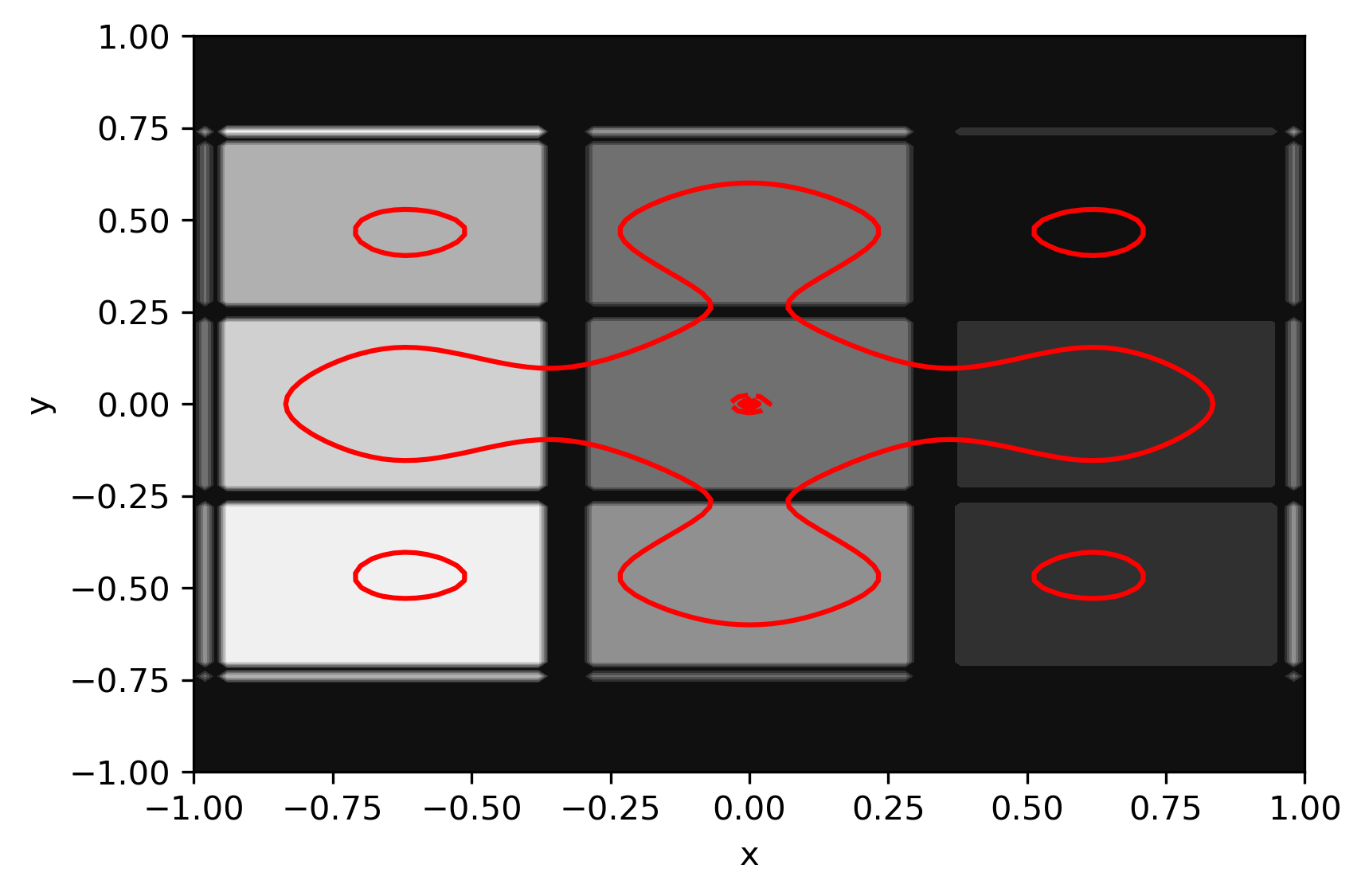} \\

        \raisebox{7em}{\rotatebox[origin=c]{90}{\textbf{With Regularization}}} &
        \includegraphics[width=0.45\textwidth]{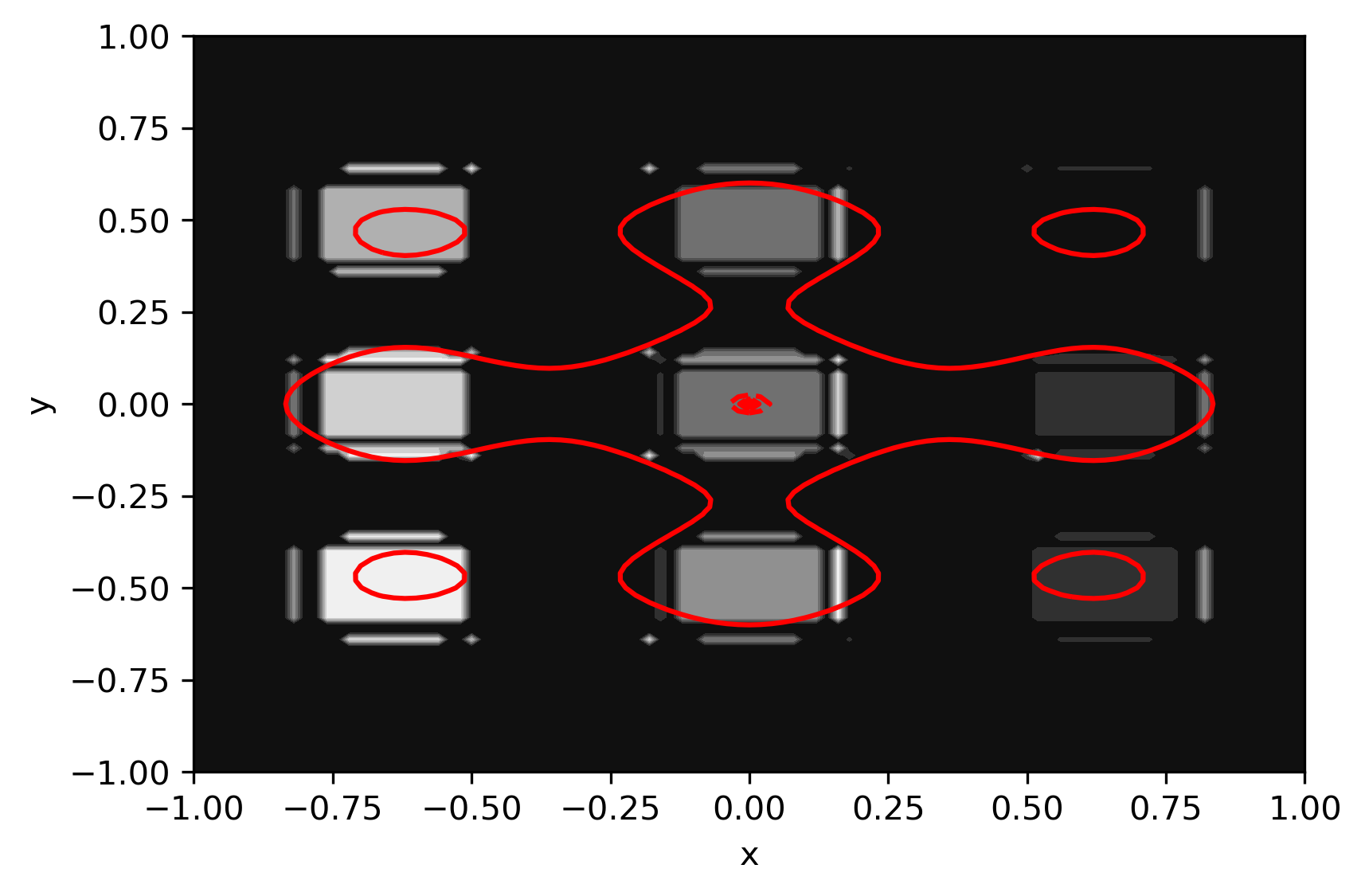} &
        \includegraphics[width=0.45\textwidth]{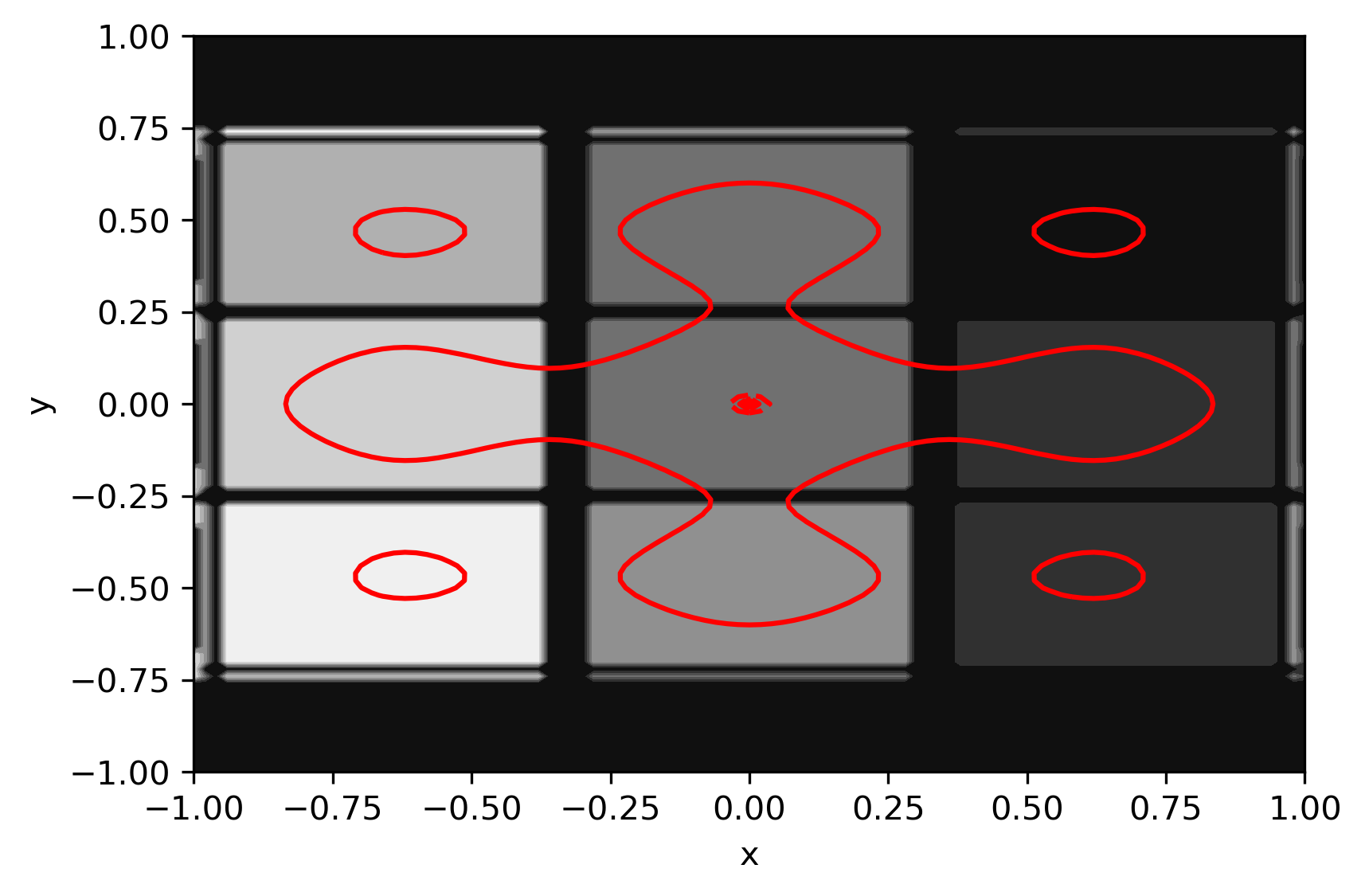} \\
    \end{tabular}

    \caption{
    Comparison of basins of attraction for second- and third-order methods for the Bohachevsky function, with and without regularization. Black regions indicate initializations that did not converge within 100 iterations. Otherwise, differently colored regions indicate initial starting points which converged to the same point.
    }
    \label{Fig: Bohachevsky Fractals}
\end{figure}

The Bohachevsky function is characterized by having many regularly spaced local minima, which is reflected in the fractals for the Bohachevsky function having regularly spaced rectangular regions. Relatively speaking, the basins of attraction for each of the local minima are significantly larger for the third-order method, with smaller regions in between. The second-order Newton algorithm has a tendency to converge to the local maxima between the local minima for initializations close to them, as it does not distinguish between different types of critical points. Because there are so many local minima and maxima, the regularization does not seem to have a discernible effect on the convergence.

\begin{figure}[H]
    \centering

    \begin{tabular}{c c c}
        & \textbf{Second-order} & \textbf{Third-order} \\

        \raisebox{7em}{\rotatebox[origin=c]{90}{\textbf{Without Regularization}}} &
        \includegraphics[width=0.45\textwidth]{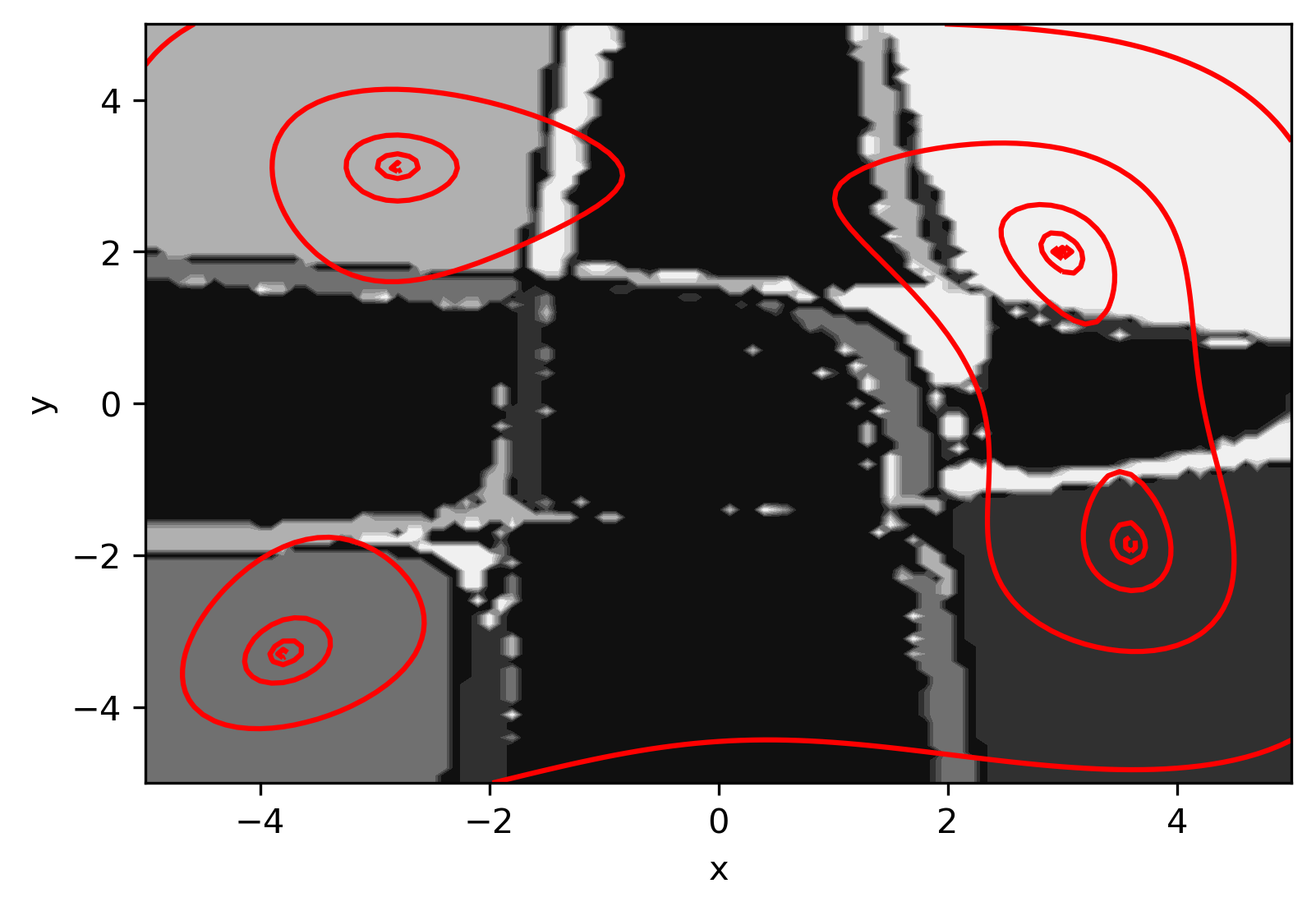} &
        \includegraphics[width=0.45\textwidth]{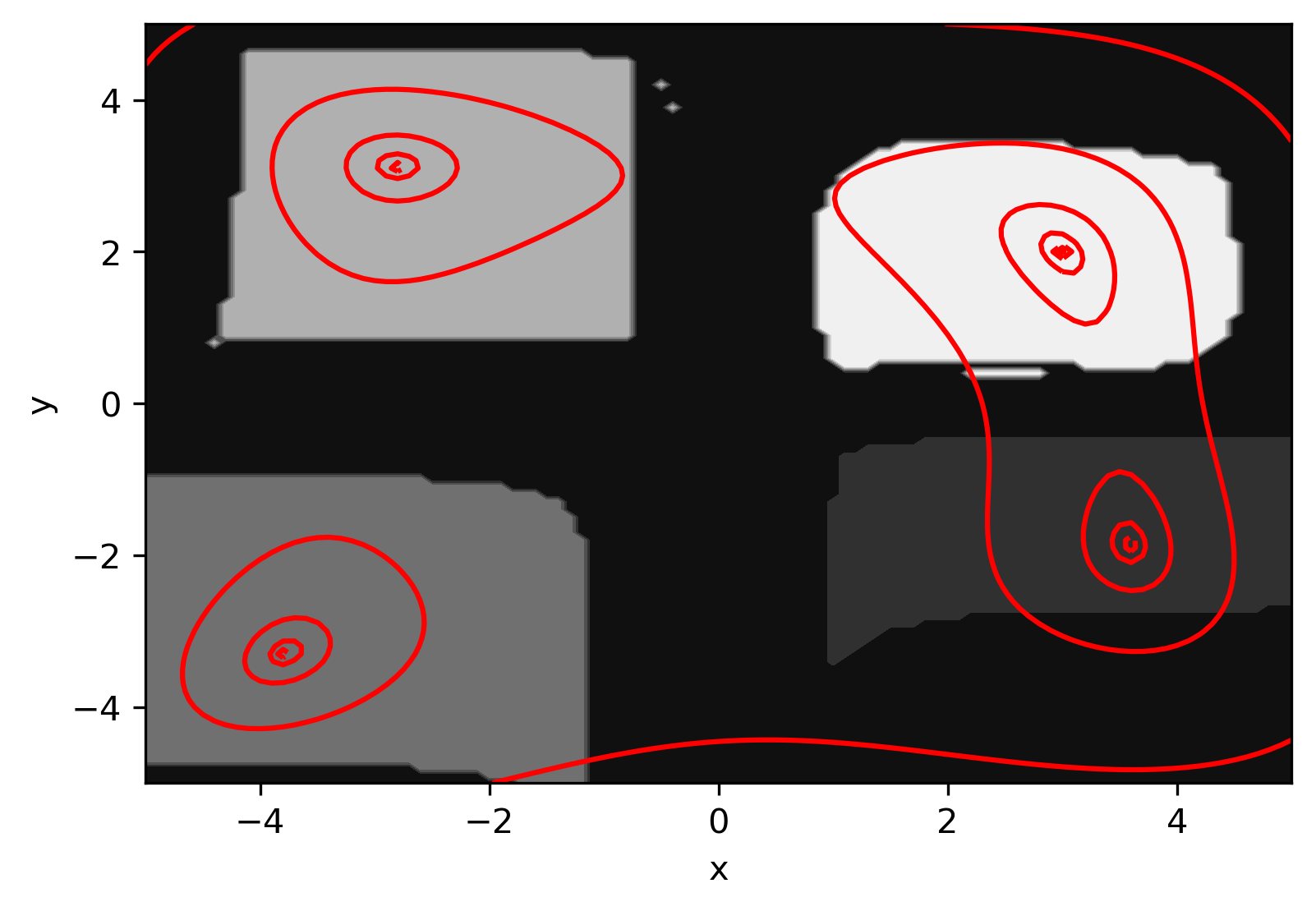} \\

        \raisebox{7em}{\rotatebox[origin=c]{90}{\textbf{With Regularization}}} &
        \includegraphics[width=0.45\textwidth]{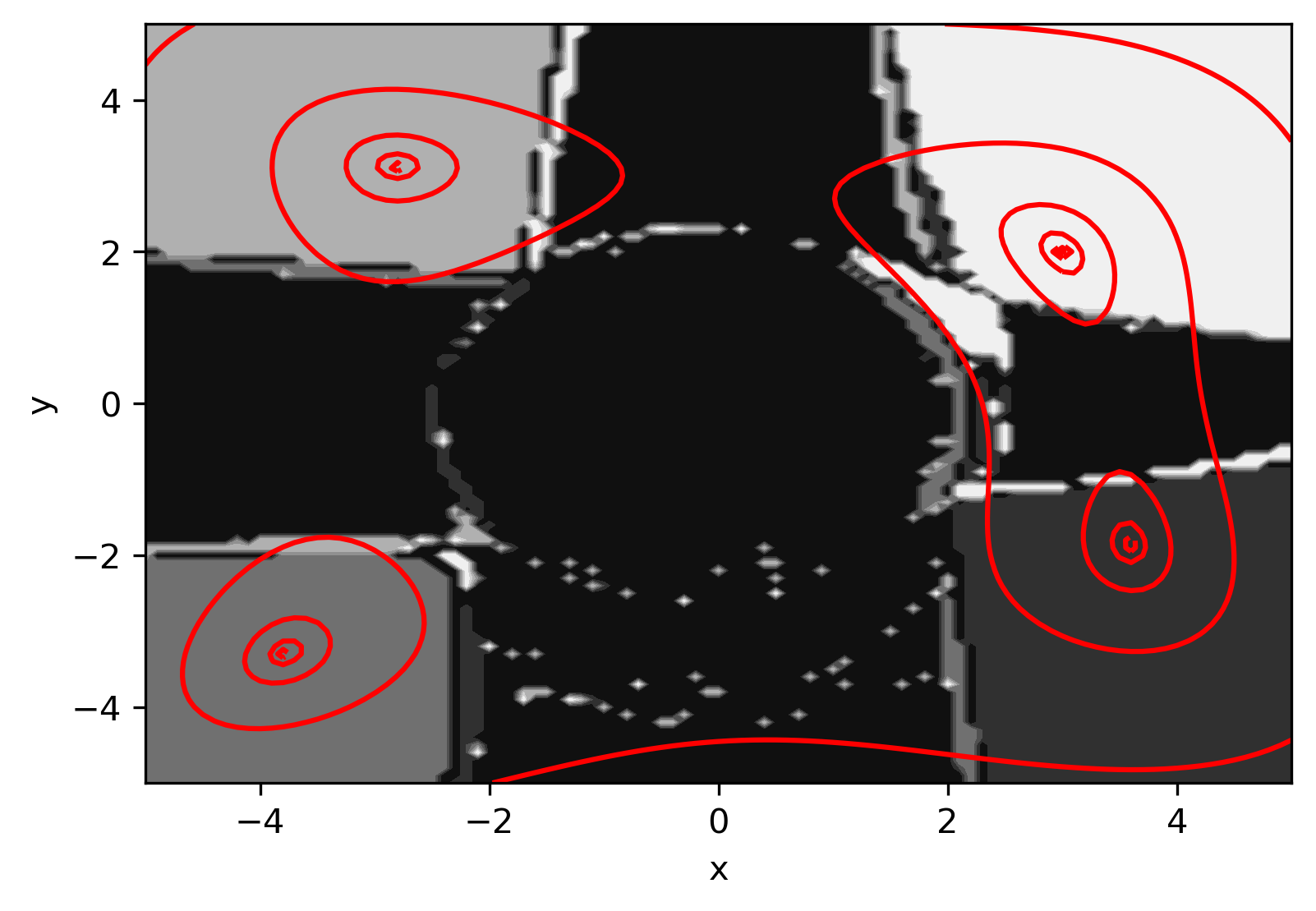} &
        \includegraphics[width=0.45\textwidth]{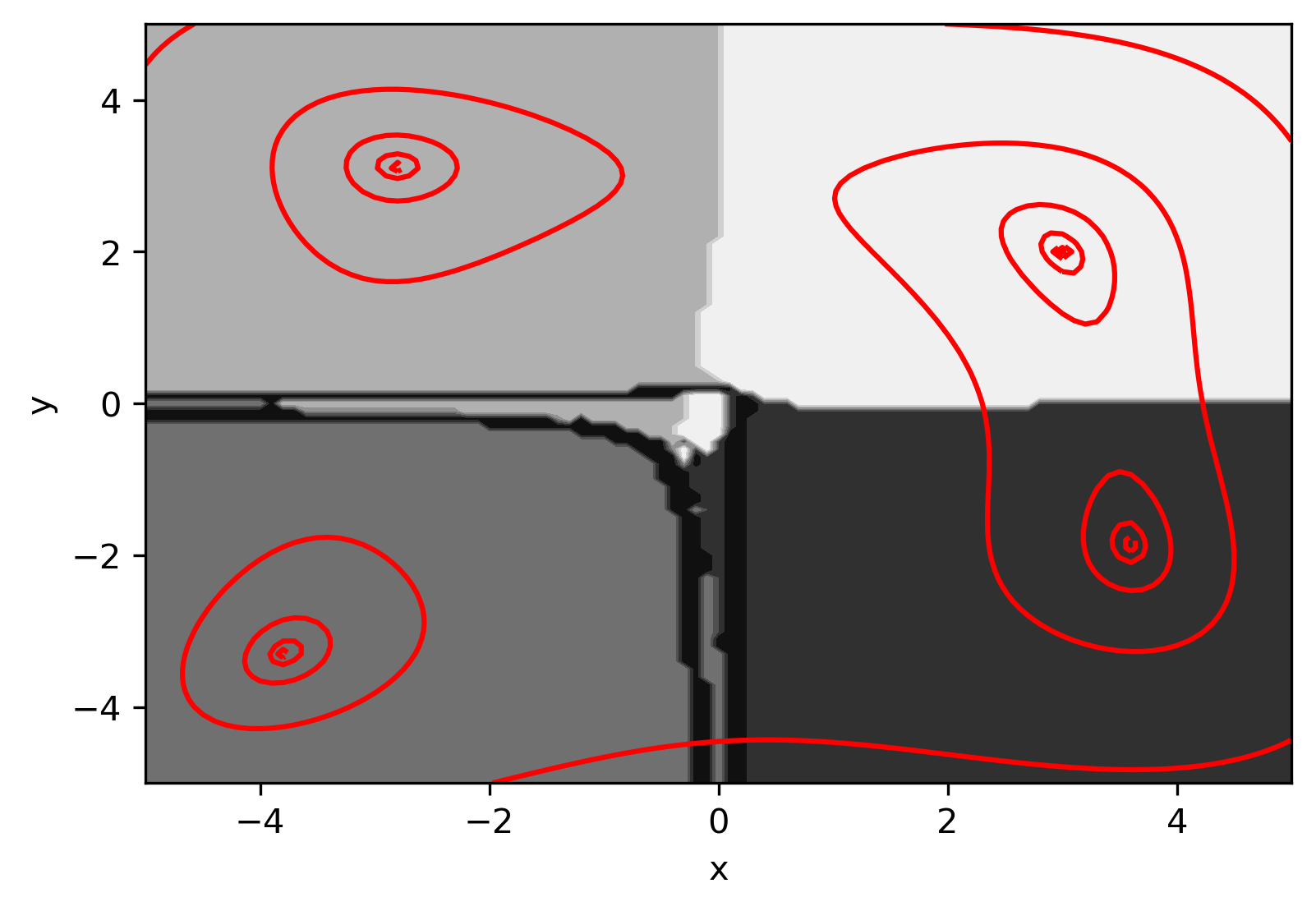} \\
    \end{tabular}

    \caption{
    Comparison of basins of attraction for second- and third-order methods for the Himmelblau function, with and without regularization. Black regions indicate initializations that did not converge within 100 iterations. Otherwise, differently colored regions indicate initial starting points which converged to the same point.
    }
    \label{Fig: Himmelblau Fractals}
\end{figure}

The Himmelblau fractals demonstrate the most difference between regularized and unregularized fractals. The most difficult region for the methods was the cross-shaped region on the $x$ and $y$ axes. Compared to the second-order method, the third-order method without regularization had narrower failure regions between the local minima, but had smaller basins of convergence. However, once regularization was added to the third-order method, the failure region disappeared almost entirely.

\begin{figure}[H]
    \centering

    \begin{tabular}{c c c}
        & \textbf{Second-order} & \textbf{Third-order} \\

        \raisebox{7em}{\rotatebox[origin=c]{90}{\textbf{Without Regularization}}} &
        \includegraphics[width=0.45\textwidth]{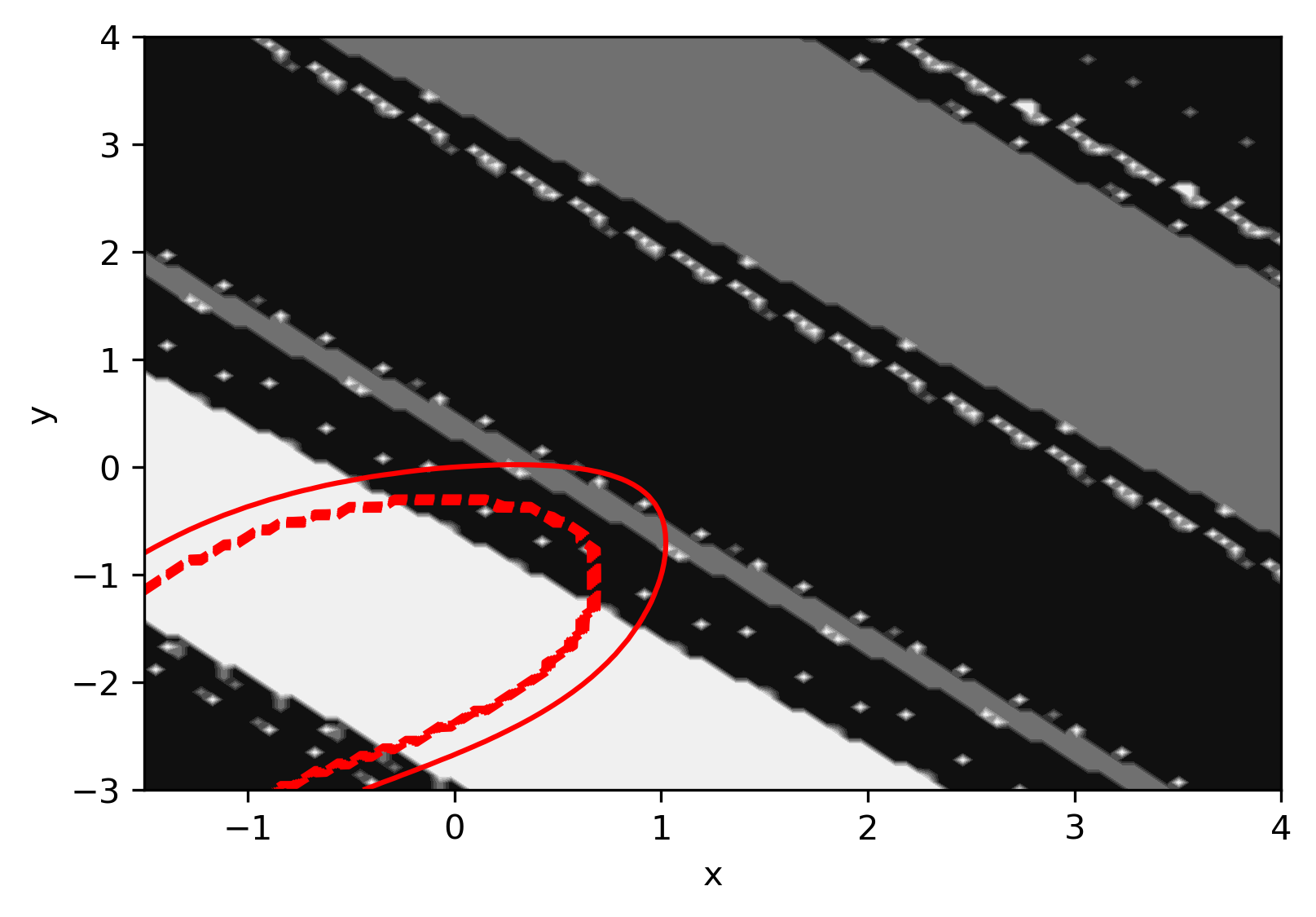} &
        \includegraphics[width=0.45\textwidth]{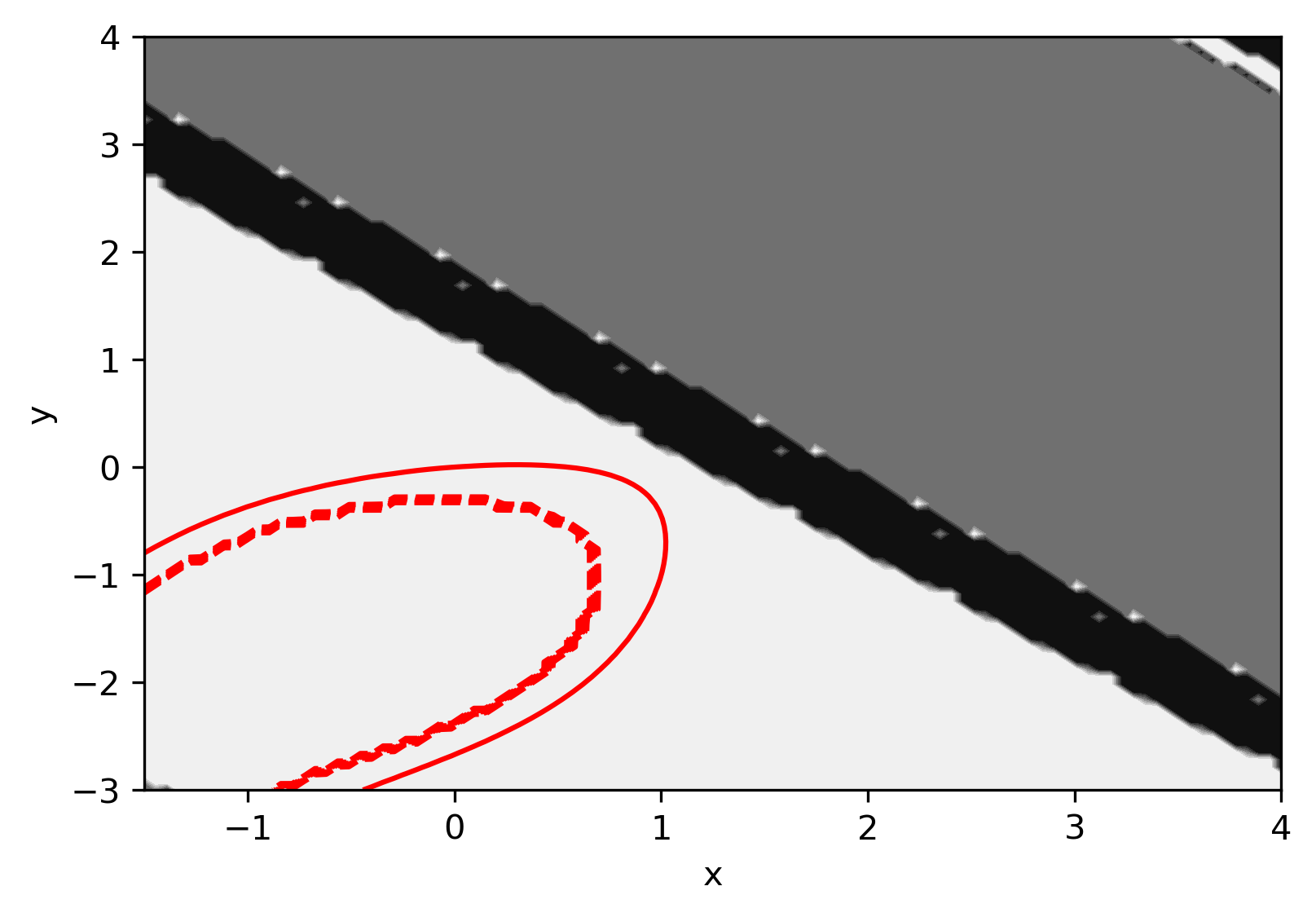} \\

        \raisebox{7em}{\rotatebox[origin=c]{90}{\textbf{With Regularization}}} &
        \includegraphics[width=0.45\textwidth]{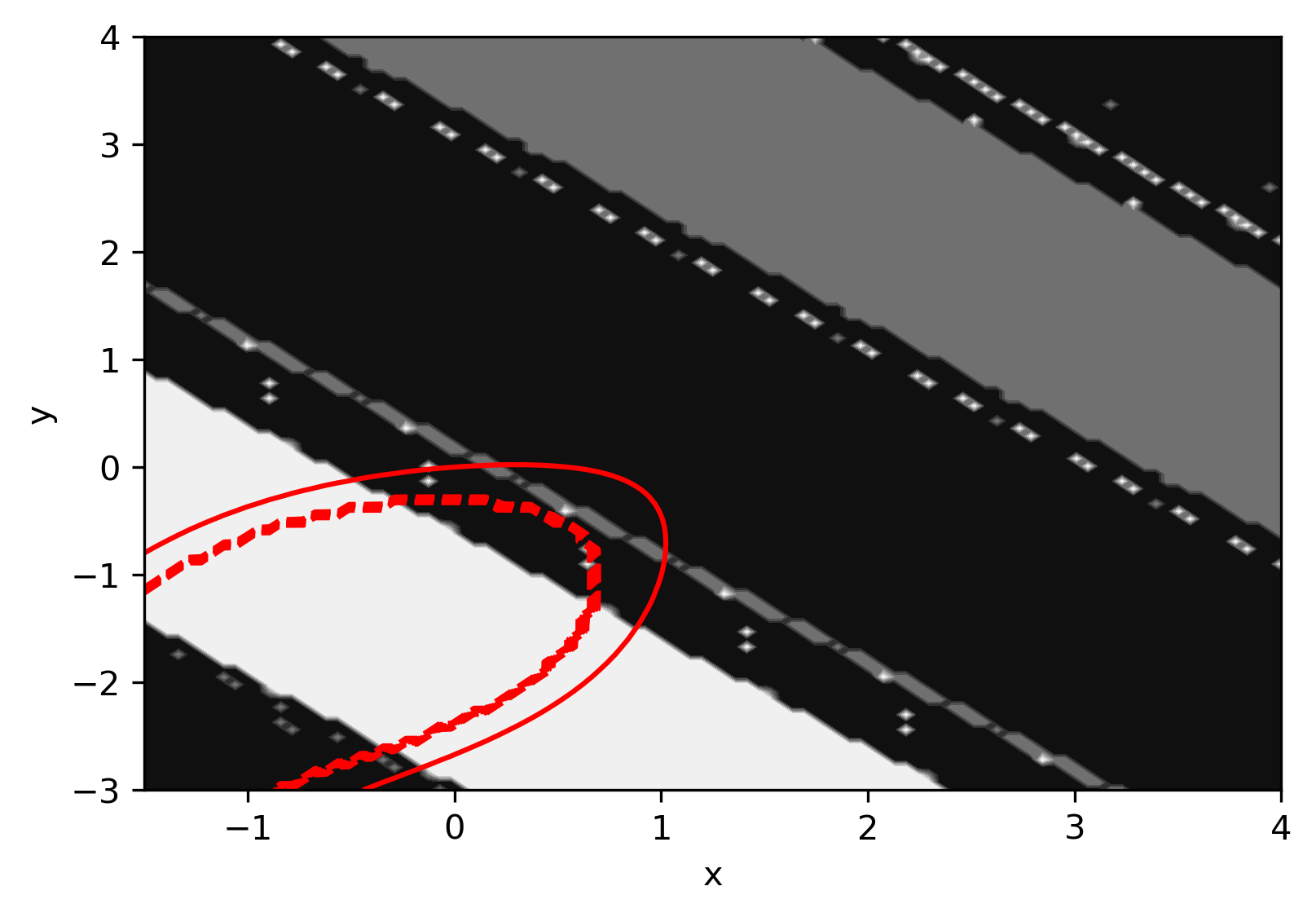} &
        \includegraphics[width=0.45\textwidth]{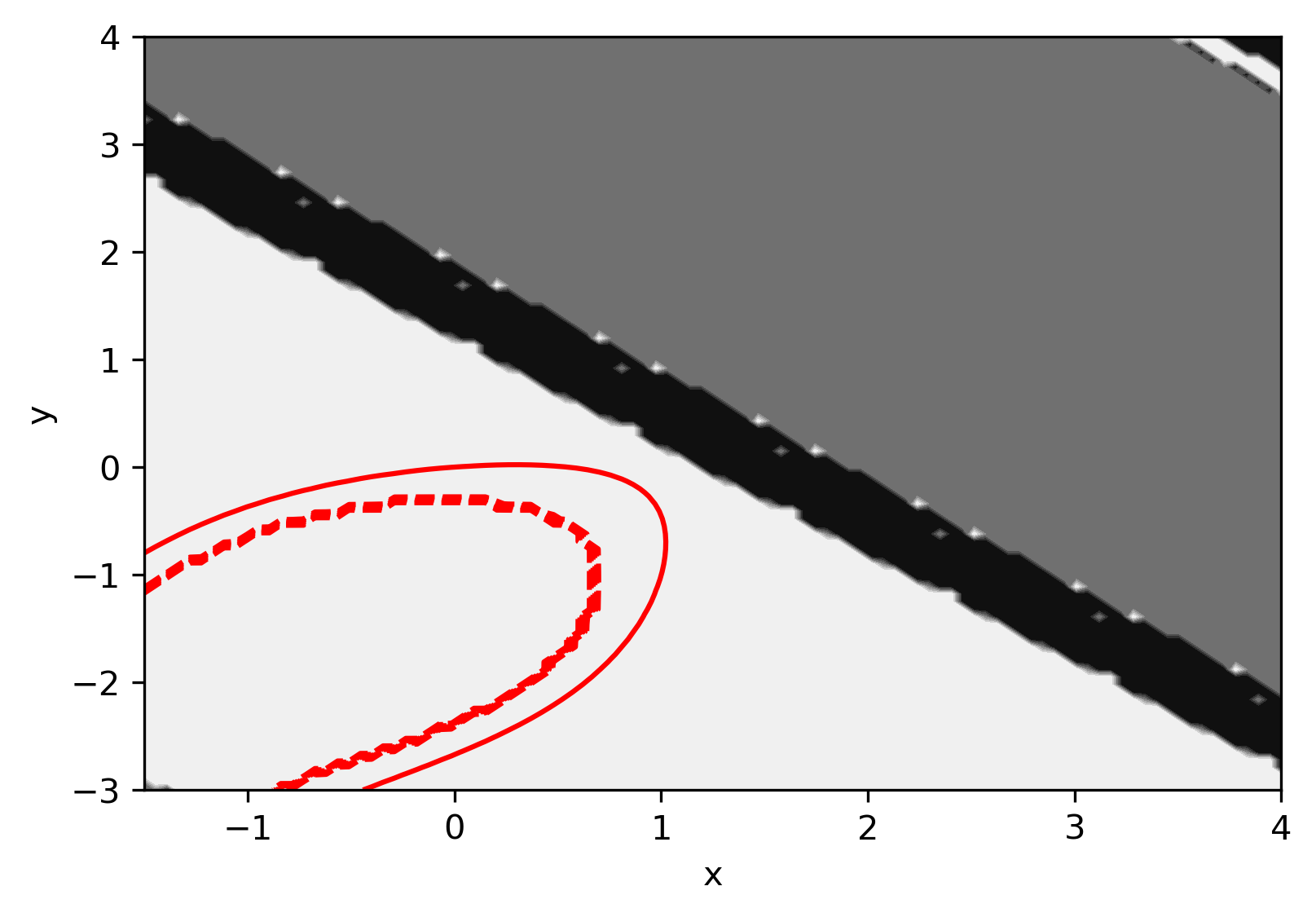} \\
    \end{tabular}

    \caption{
    Comparison of basins of attraction for second- and third-order methods for the McCormick function, with and without regularization. Black regions indicate initializations that did not converge within 100 iterations. Otherwise, differently colored regions indicate initial starting points which converged to the same point.
    }
    \label{Fig: McCormock Fractals}
\end{figure}

Similar to the Bohachevsky function, the McCormick function can be characterized by regularly spaced local minima. We observe similar behavior as well; for both the second and third-order methods the fractals can be characterized as bands that correspond to a local mininmum, with little effect of regularization. We again observe that the basins are much larger for the third order methods, and again the second-order method features regions that converge to a local maximum rather than a local minimum.

\section{Conclusion}\label{Sec: Conclusion}

In this work, we developed an unregularized third-order Newton method and established theoretical guarantees for its behavior, supported by numerical experiments. Compared to the classical Newton method, the primary advantage lies in its wider basins of attraction and reduced sensitivity to initialization. However, from a computational standpoint, the method is not competitive with second-order approaches in regimes where those second-order methods converge reliably. The dominant cost arises from solving a semidefinite program at each iteration, and the resulting reduction in iteration count does not generally offset this overhead. With that said, integrating third-order information is certainly worth consideration in settings where second-order methods may not be reliable.

These limitations suggest several directions for future work. One natural approach is to relax the requirement of solving each SDP to high precision. It is plausible that approximate solutions—obtained, for example, by terminating the SDP solver early—are sufficient to preserve the qualitative behavior of the method while significantly reducing computational cost. Alternatively, one can consider interior-point type iterations for the formulation posed in (\ref{Eq: minimize cubic}), given that it is in principle a convex problem with a well-established barrier function. A more ambitious direction is to develop third-order methods that retain the favorable geometric properties observed here without requiring an SDP solve. For instance, one could consider approximations that use only partial third-order information, such as diagonal or structured subsets of the third-derivative tensor, in the spirit of quasi-Newton methods. This raises several questions: whether local minima of the resulting models can be computed efficiently, whether convergence guarantees can be maintained, and whether the enlarged basins of attraction persist. One can also consider integrating advances in scalable semidefinite programming~\cite{majumdar2020recent}, though theoretically minimizing cubic polynomials is in general as hard as SDP~\cite{ahmadi2022complexity}.

More broadly, this work suggests a framework for a class of third-order methods that treat minimizing cubic polynomials as a blackbox. We again refer the reader to the work in \cite{cai2026globally} that globalizes this framework. One can also study classic adaptations of second-order methods beyond Levenberg–Marquardt regularization which was studied in this work, such as trust regions and exact line search. We do point out that minimizing a cubic polynomial over a sphere is NP-hard~\cite{nesterov2003random}, though one can easily add a norm constraint into the SDPs in this work.
\section*{Conflict of interest}
The authors declare that they have no conflict of interest.

\section*{Data availability}
All code and data are available publicly at \url{https://github.com/jeffreyzhang92/Third_Order_Newton}.

\bibliographystyle{abbrv}
\bibliography{refs}

\begin{appendix}
    \section{Complete Newton Fractals}\label{Sec: Complete Newton Fractals}
\end{appendix}

\subsection{Beale Function}
\begin{figure}[H]
    \centering

    \begin{tabular}{c c c}
        & \textbf{Cubic Objective} & \textbf{Gap Objective} \\

        \raisebox{7em}{\rotatebox[origin=c]{90}{\textbf{LM Regularization}}} &
        \includegraphics[width=0.45\textwidth]{imgs/Beale_3_LM_cubic_Fractal.png} &
        \includegraphics[width=0.45\textwidth]{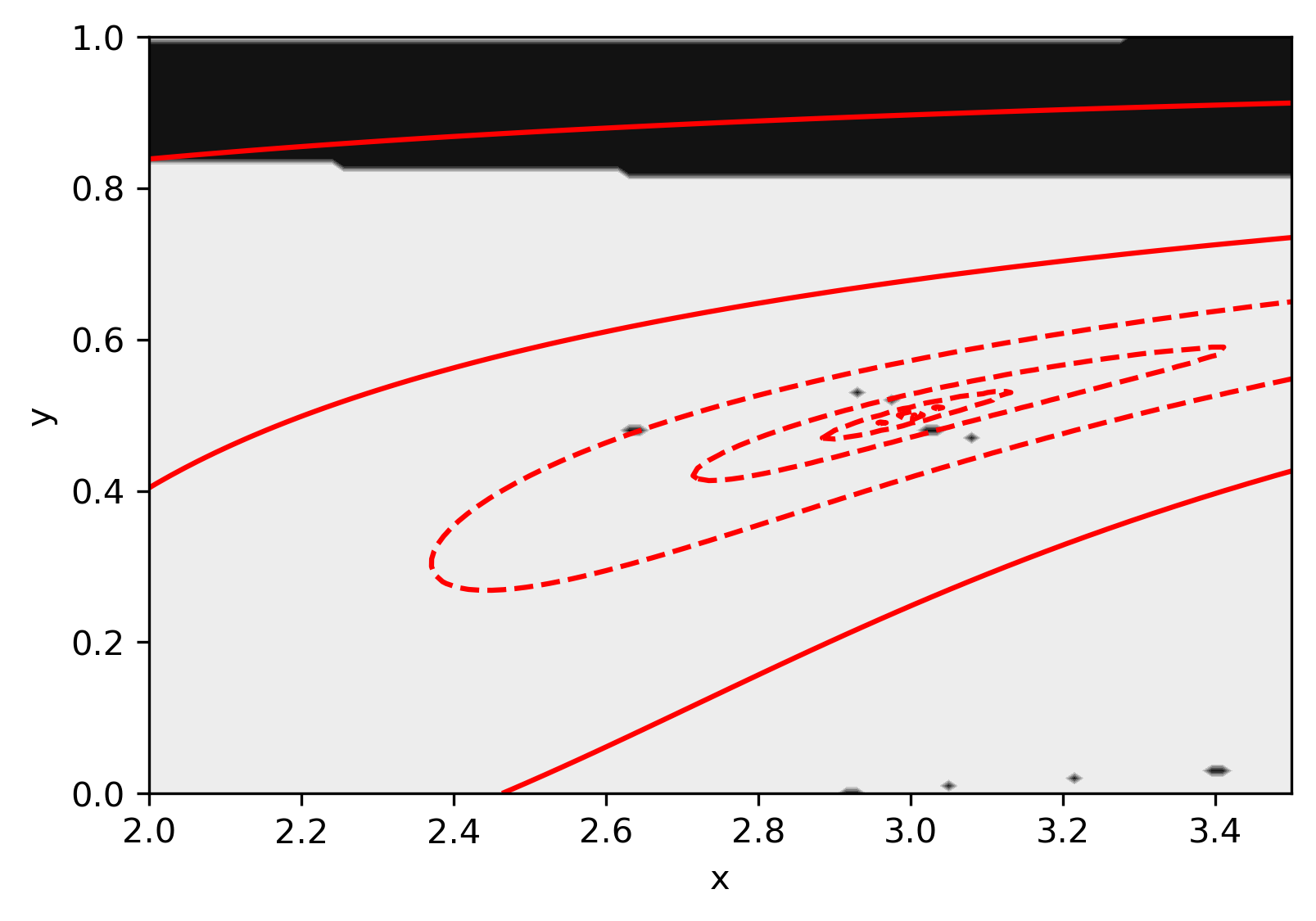} \\

        \raisebox{7em}{\rotatebox[origin=c]{90}{\textbf{Heuristic Regularization}}} &
        \includegraphics[width=0.45\textwidth]{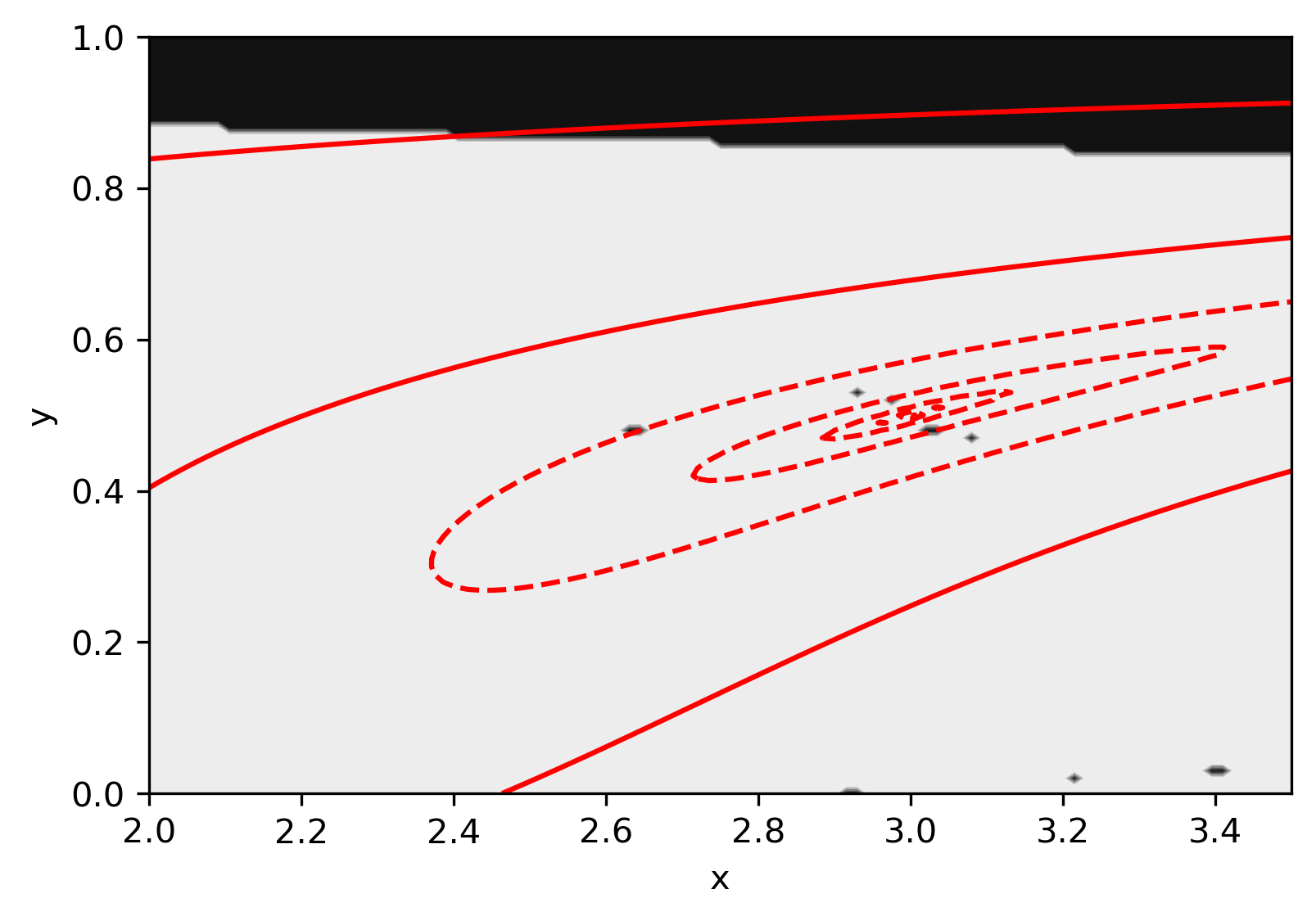} &
        \includegraphics[width=0.45\textwidth]{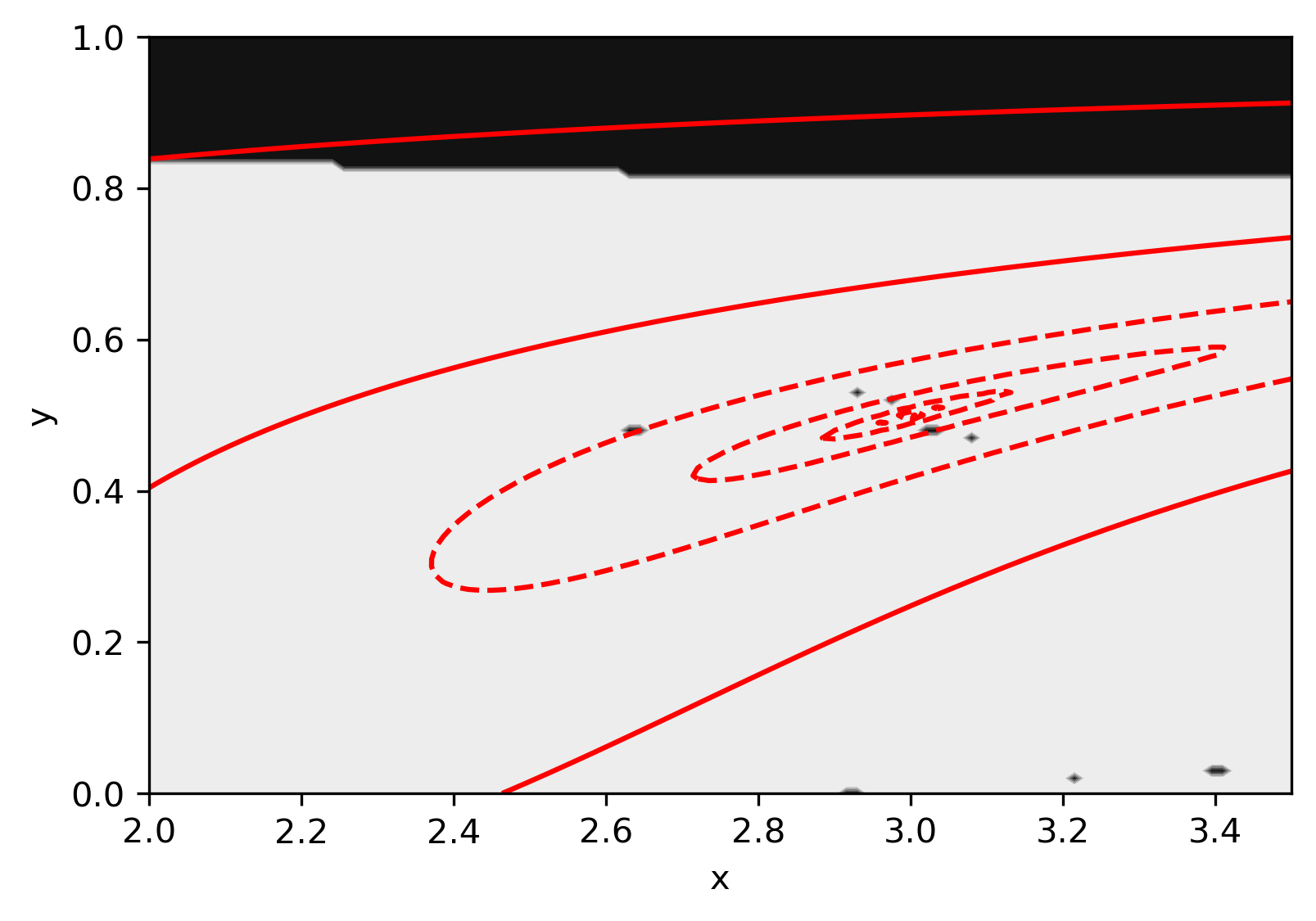} \\
    \end{tabular}

    \caption{
    Comparison of basins of attraction for different specifications of the third-order method for the Bohachevsky function, with and without regularization. Black regions indicate initializations that did not converge within 100 iterations. Otherwise, differently colored regions indicate initial starting points which converged to the same point.
    }
    \label{Fig: Beale Fractals 3}
\end{figure}

\subsection{Bohachevsky Function}

\begin{figure}[H]
    \centering

    \begin{tabular}{c c c}
        & \textbf{Cubic Objective} & \textbf{Gap Objective} \\

        \raisebox{7em}{\rotatebox[origin=c]{90}{\textbf{LM Regularization}}} &
        \includegraphics[width=0.45\textwidth]{imgs/Bohachevsky_3_LM_cubic_Fractal.png} &
        \includegraphics[width=0.45\textwidth]{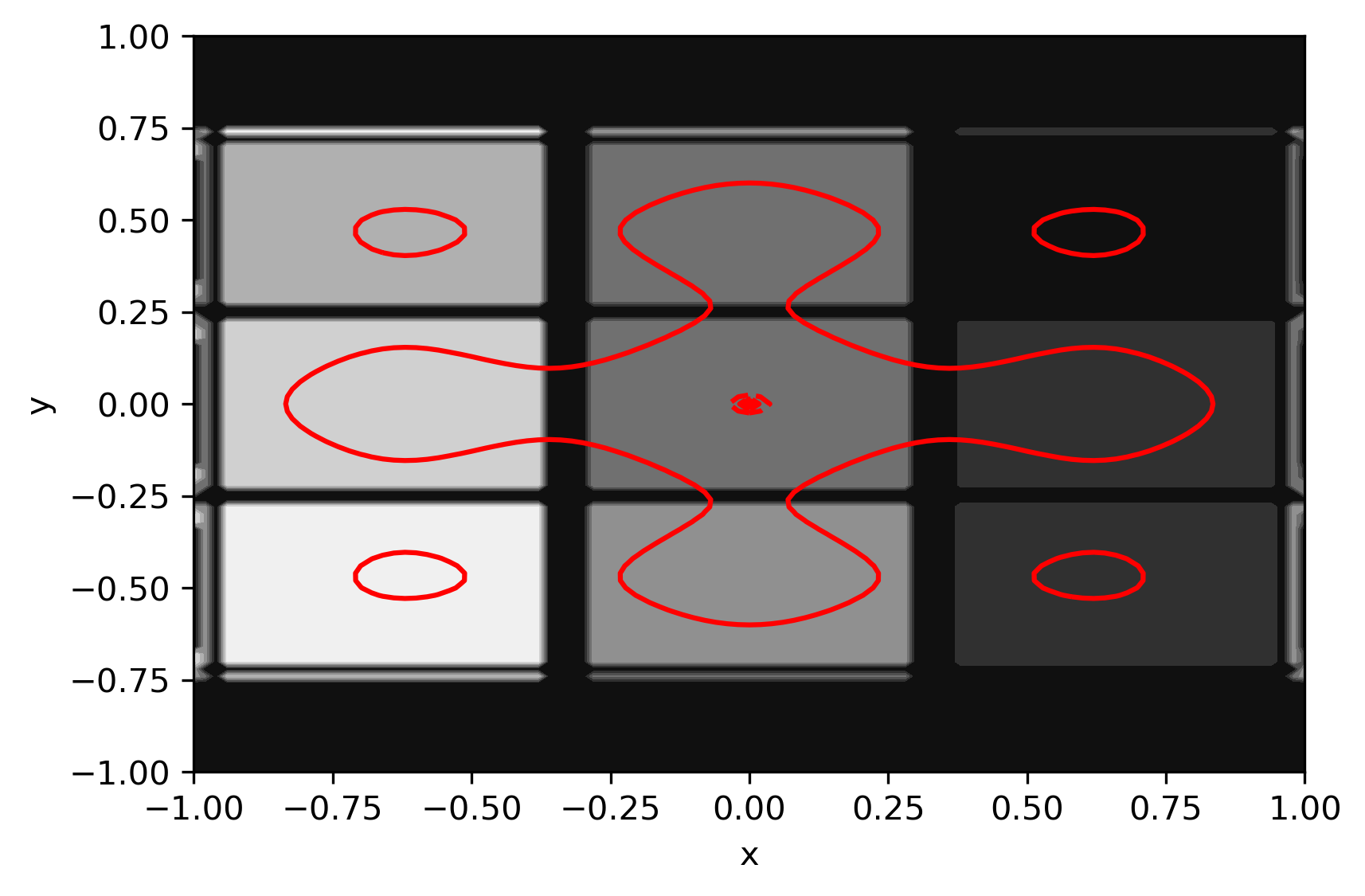} \\

        \raisebox{7em}{\rotatebox[origin=c]{90}{\textbf{Heuristic Regularization}}} &
        \includegraphics[width=0.45\textwidth]{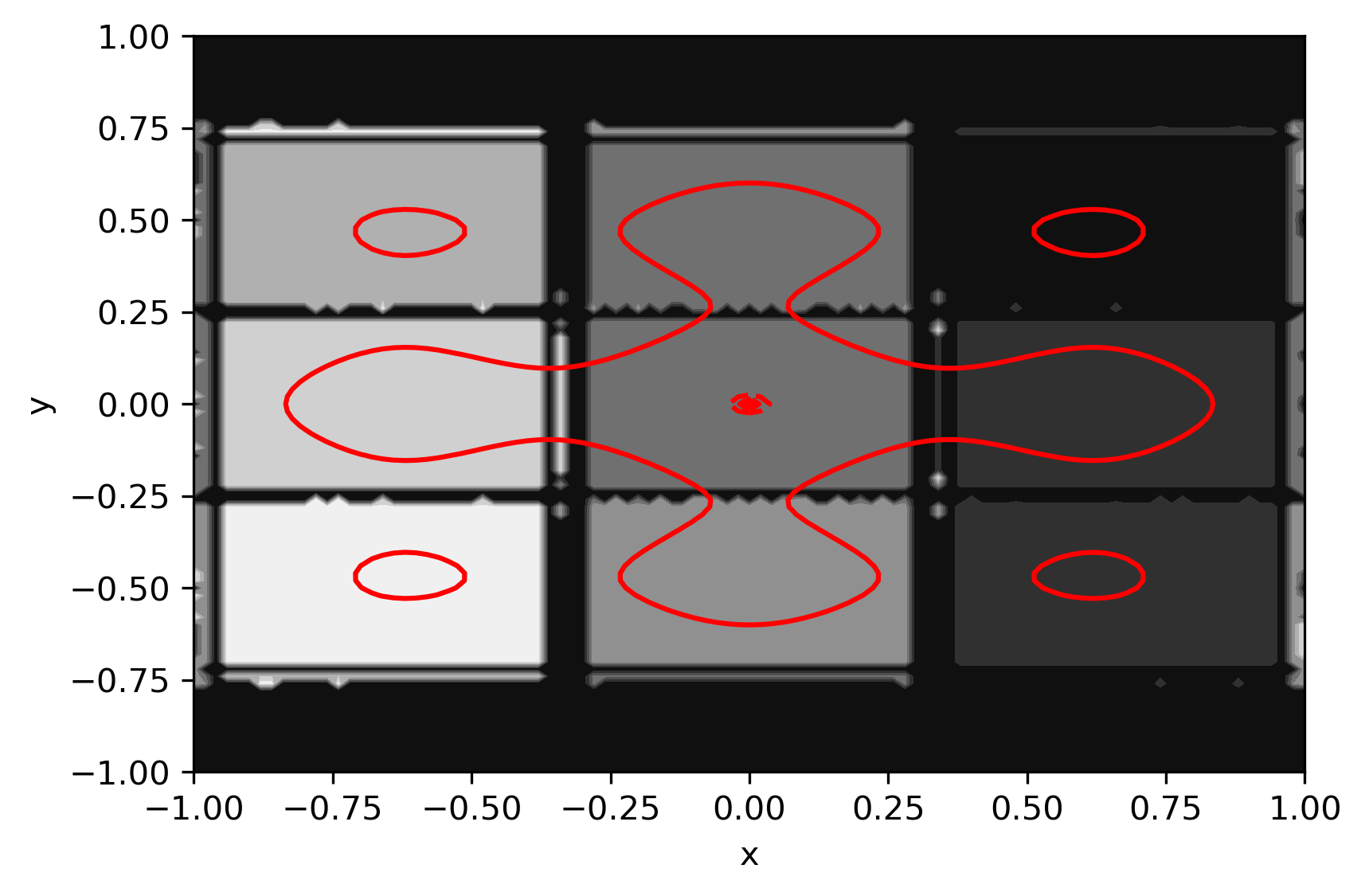} &
        \includegraphics[width=0.45\textwidth]{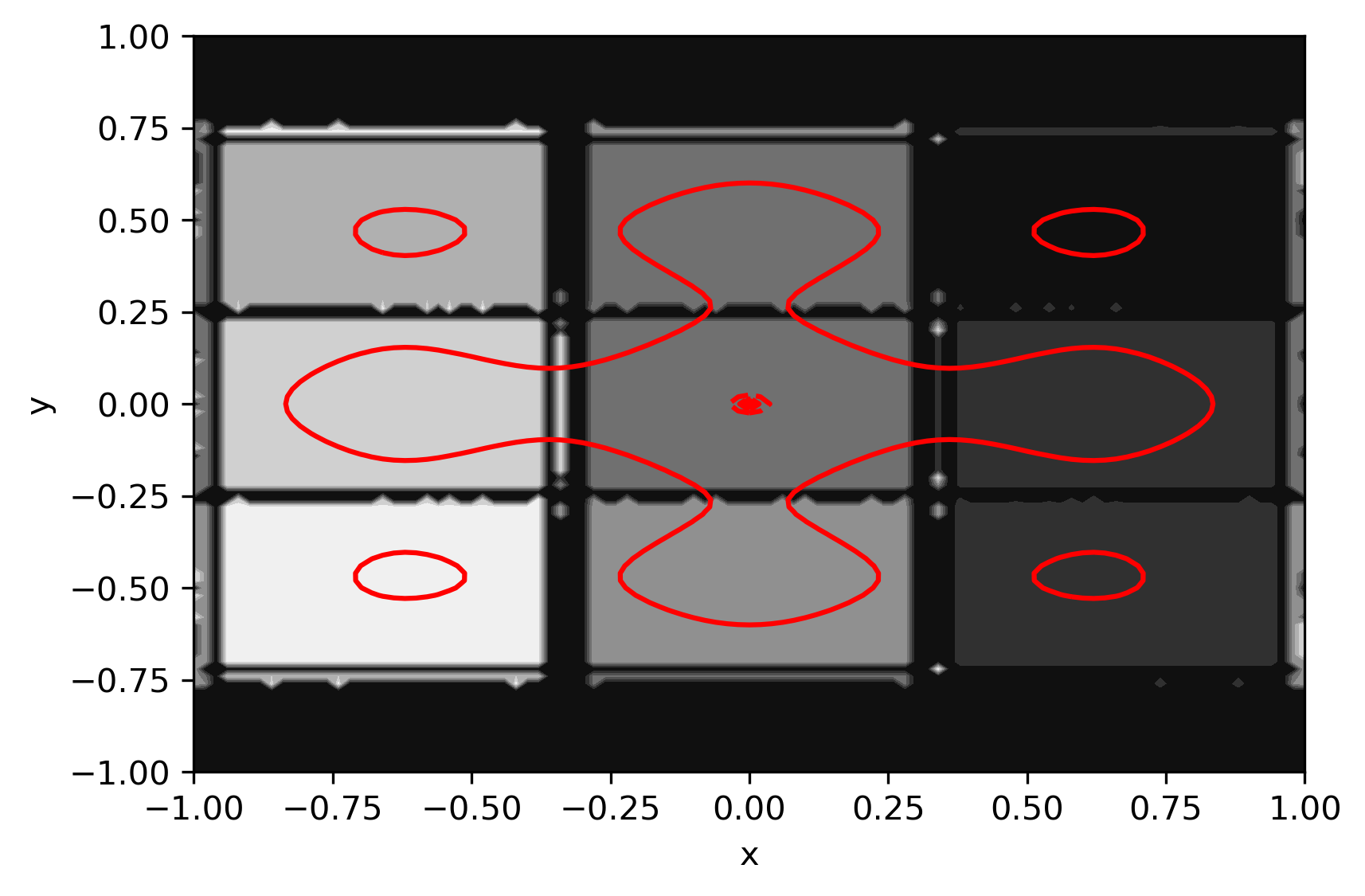} \\
    \end{tabular}

    \caption{
    Comparison of basins of attraction for different specifications of the third-order method for the Bohachevsky function, with and without regularization. Black regions indicate initializations that did not converge within 100 iterations. Otherwise, differently colored regions indicate initial starting points which converged to the same point.
    }
    \label{Fig: Bohachevsky Fractals 3}
\end{figure}

\subsection{Himmelblau Function}

\begin{figure}[H]
    \centering

    \begin{tabular}{c c c}
        & \textbf{Cubic Objective} & \textbf{Gap Objective} \\

        \raisebox{7em}{\rotatebox[origin=c]{90}{\textbf{LM Regularization}}} &
        \includegraphics[width=0.45\textwidth]{imgs/Himmelblau_3_LM_cubic_Fractal.png} &
        \includegraphics[width=0.45\textwidth]{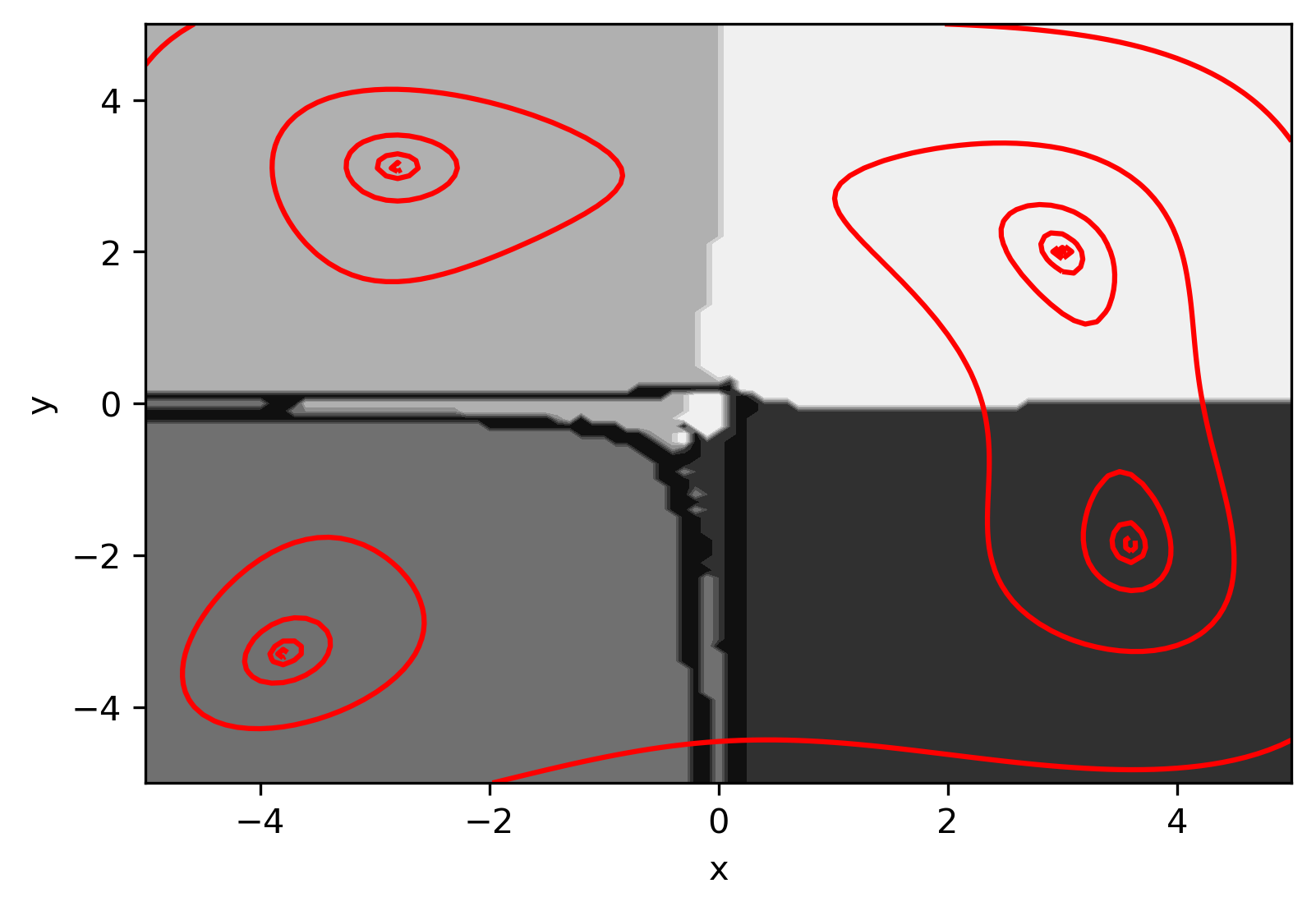} \\

        \raisebox{7em}{\rotatebox[origin=c]{90}{\textbf{Heuristic Regularization}}} &
        \includegraphics[width=0.45\textwidth]{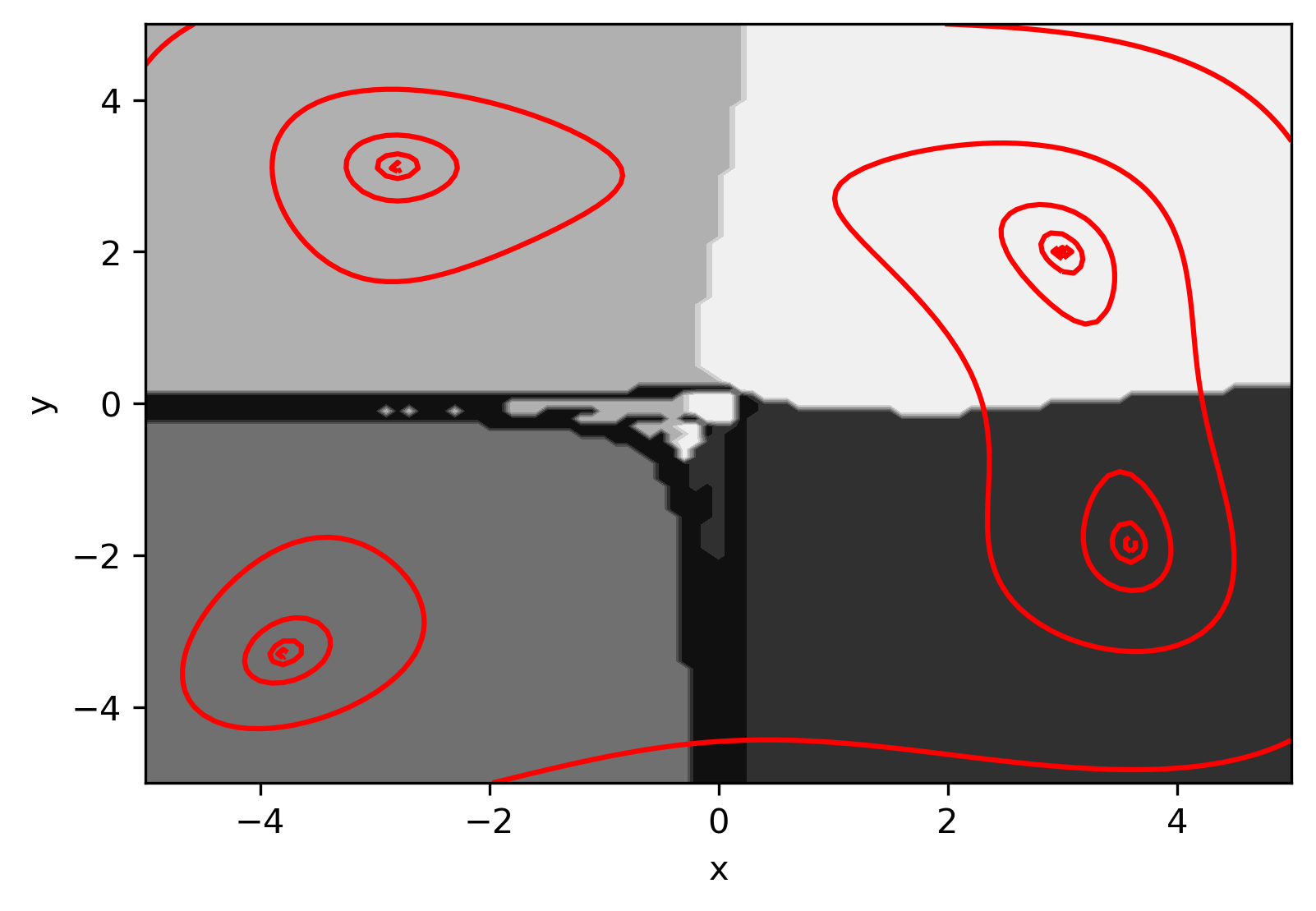} &
        \includegraphics[width=0.45\textwidth]{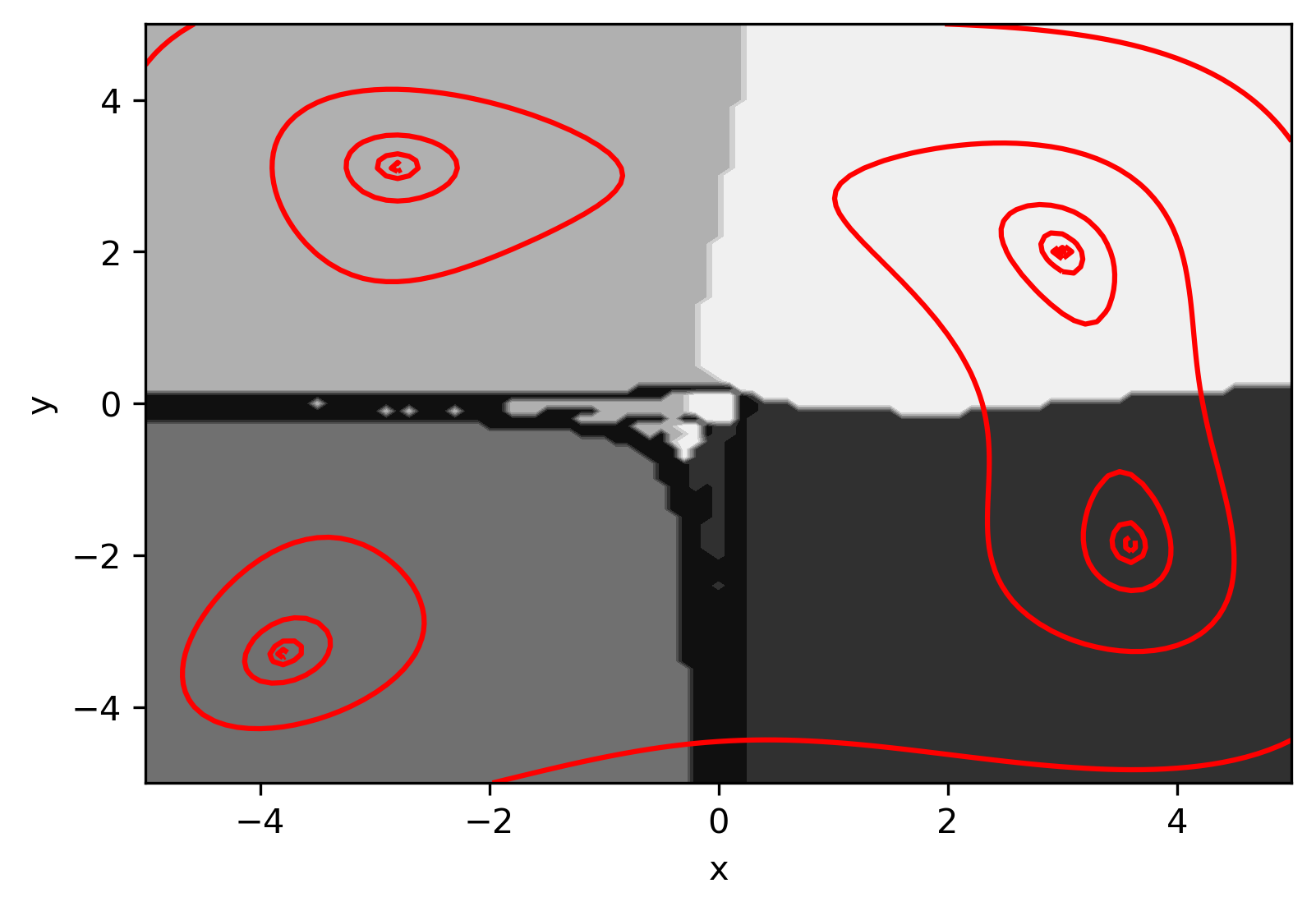} \\
    \end{tabular}

    \caption{
    Comparison of basins of attraction for different specifications of the third-order method for the Himmelblau function, with and without regularization. Black regions indicate initializations that did not converge within 100 iterations. Otherwise, differently colored regions indicate initial starting points which converged to the same point.
    }
    \label{Fig: Himmelblau Fractals 3}
\end{figure}

\subsection{McCormick Function}

\begin{figure}[H]
    \centering

    \begin{tabular}{c c c}
        & \textbf{Cubic Objective} & \textbf{Gap Objective} \\

        \raisebox{7em}{\rotatebox[origin=c]{90}{\textbf{LM Regularization}}} &
        \includegraphics[width=0.45\textwidth]{imgs/McCormick_3_LM_cubic_Fractal.png} &
        \includegraphics[width=0.45\textwidth]{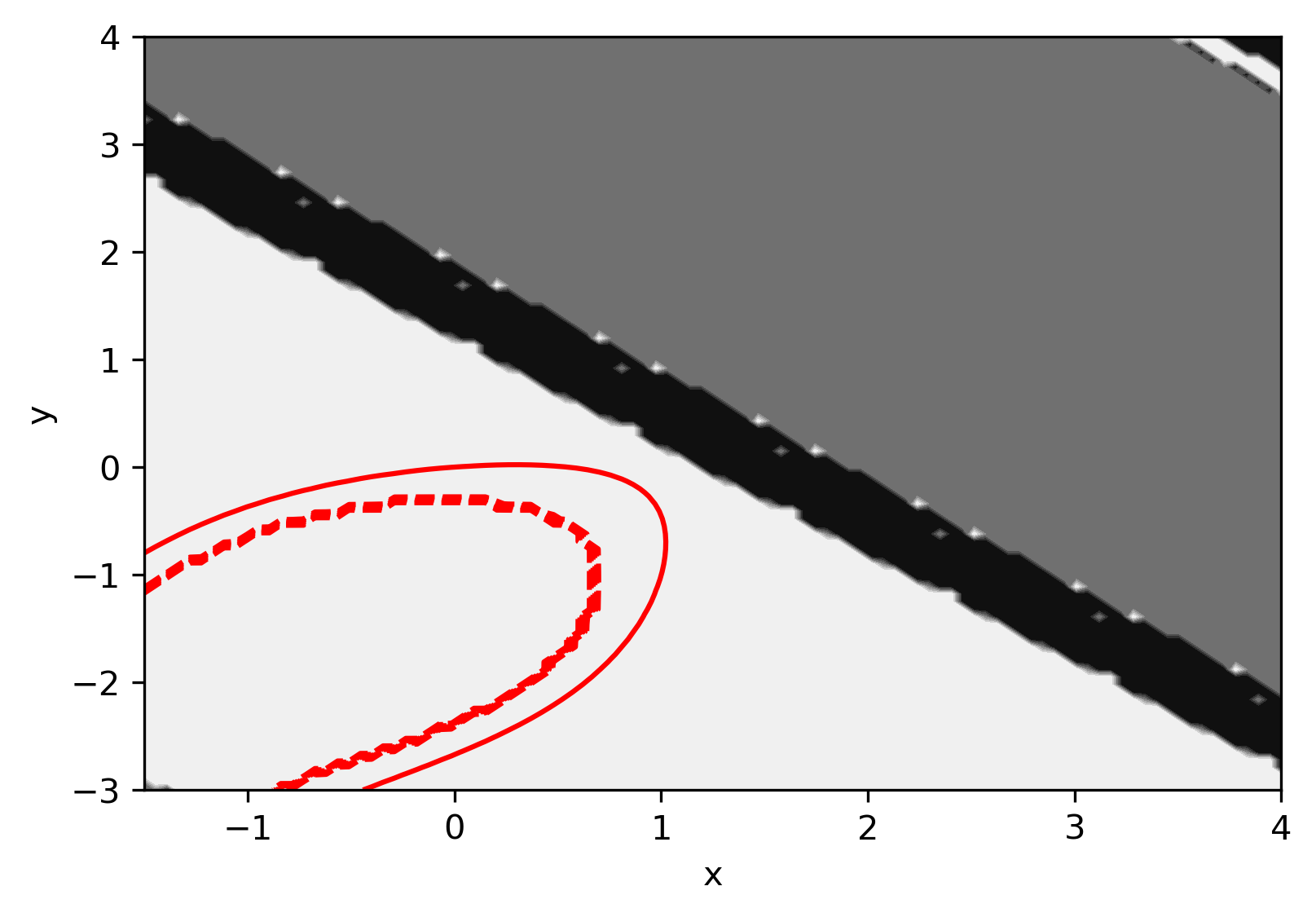} \\

        \raisebox{7em}{\rotatebox[origin=c]{90}{\textbf{Heuristic Regularization}}} &
        \includegraphics[width=0.45\textwidth]{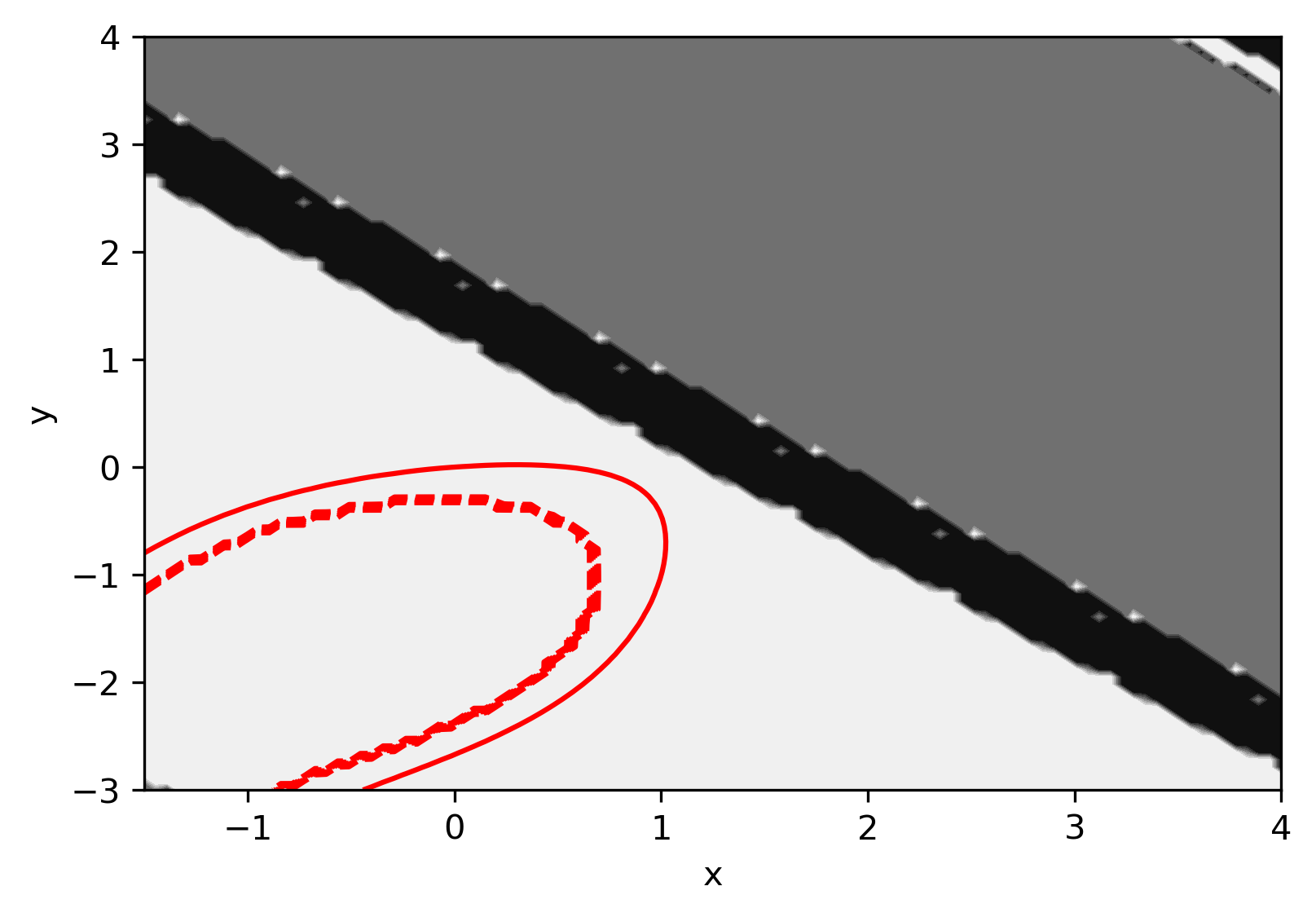} &
        \includegraphics[width=0.45\textwidth]{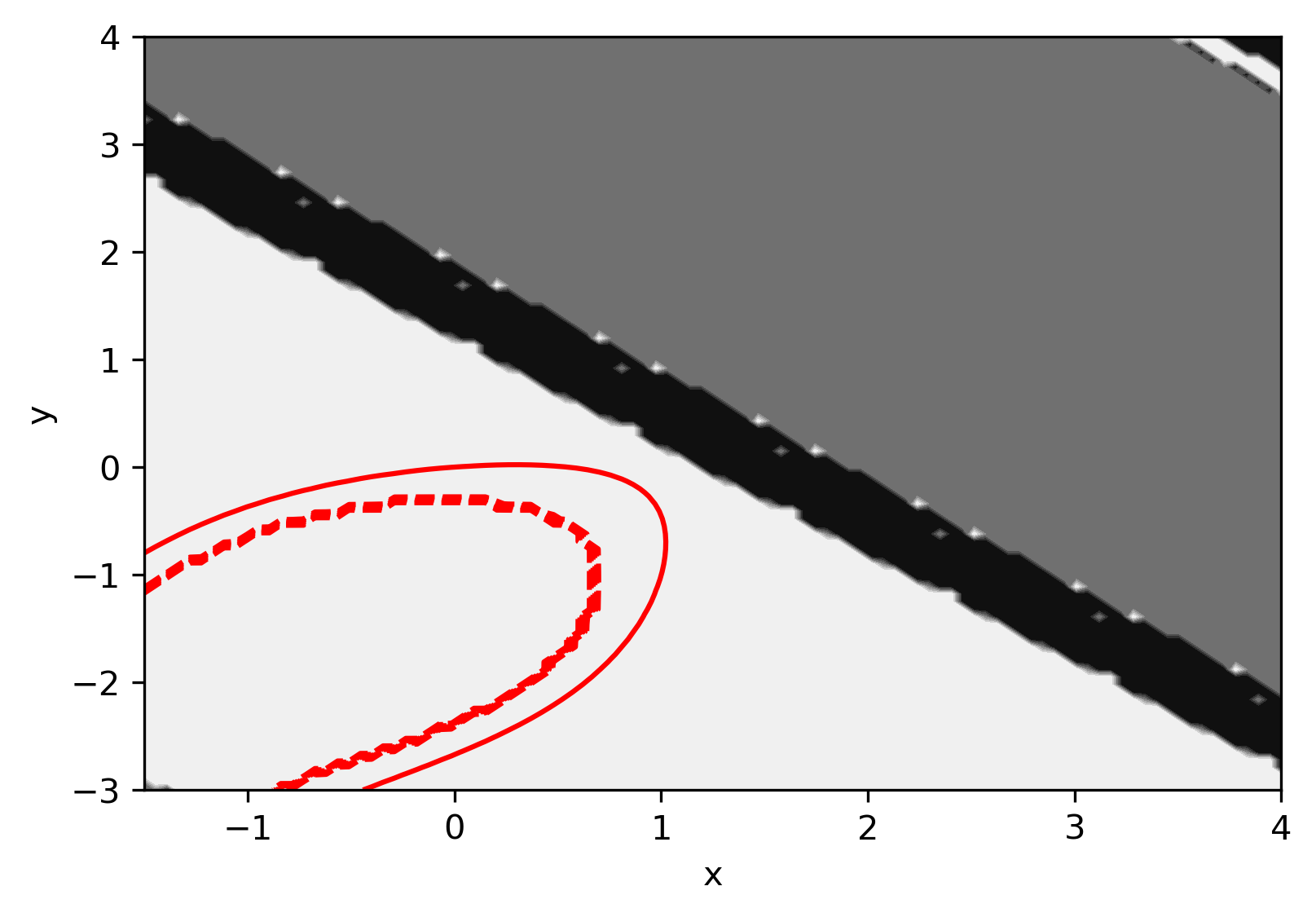} \\
    \end{tabular}

    \caption{
    Comparison of basins of attraction for different specifications of the third-order method for the McCormick function, with and without regularization. Black regions indicate initializations that did not converge within 100 iterations. Otherwise, differently colored regions indicate initial starting points which converged to the same point.
    }
    \label{Fig: McCormick Fractals 3}
\end{figure}

\end{document}